\documentclass[11pt]{article}
\usepackage[margin=1in]{geometry}
\usepackage[T1]{fontenc}
\usepackage{lmodern}
\usepackage{amsmath,amssymb,amsthm,bm}
\usepackage{mathrsfs}
\usepackage{algorithm,algpseudocode}
\usepackage{booktabs,microtype}
\usepackage{xcolor}
\usepackage{graphicx}
\usepackage[hidelinks]{hyperref}
\hypersetup{
  pdftitle={Fast HOSVD via Randomized Preconditioning and Mixed Precision},
  pdfauthor={Yiying Chen and Yuning Yang}
}
\floatname{algorithm}{Algorithm}
\algrenewcommand\algorithmicrequire{\textbf{Input:}}
\algrenewcommand\algorithmicensure{\textbf{Output:}}
\algrenewcommand\algorithmicreturn{\textbf{return}}
\theoremstyle{definition}
\newtheorem{proposition}{Proposition}[section]
\newtheorem{lemma}[proposition]{Lemma}
\newtheorem{theorem}[proposition]{Theorem}

\newcommand{\R}{\mathbb R}
\newcommand{\range}{\operatorname{range}}

\newcommand{\diag}{\operatorname{diag}}
\newcommand{\tr}{\operatorname{tr}}

\newcommand{\norm}[1]{\left\lVert#1\right\rVert}
\newcommand{\T}{^{\mathsf T}}
\title{Fast HOSVD via Randomized Preconditioning and Mixed Precision}
\author{Yiying Chen\thanks{Faculty of English Language and Culture, Guangdong University of Foreign Studies, Guangzhou 510420, China.}
  \and Yuning Yang\thanks{Corresponding author. School of Mathematics, Guangxi University, Nanning 530004, China. Email: \mbox{\href{mailto:yyang@gxu.edu.cn}{\nolinkurl{yyang@gxu.edu.cn}}}.}}
\date{}

\begin{document}
\maketitle

\begin{abstract}
Computing dominant mode subspaces is a major cost in the higher-order singular value decomposition (HOSVD) and its sequentially truncated variant (STHOSVD). Randomized Tucker algorithms can   reduce this cost by approximating the required subspaces from small sketches, but may introduce additional approximation errors.  We exploit the extreme aspect ratios of tensor unfoldings, whose mode dimensions can remain manageable while being orders of magnitude smaller than the products of the remaining dimensions. Instead of directly approximating the desired subspaces, we sketch the entire unfolding at a scale comparable to its mode dimension to preserve its full row-space geometry, and use randomized preconditioning to recover the desired singular subspaces through small matrix factorizations. We develop randomized preconditioned HOSVD and STHOSVD algorithms incorporating numerical truncation, and further accelerate their large matrix multiplications using mixed precision. Our analysis shows that truncation adds only a term quadratic in the numerical threshold to the classical HOSVD tail-energy bound, while mixed precision introduces a relative \(O(\sqrt{u_\ell})\) correction to the tail-energy term and an additive \(O(u_\ell)\) term under suitable conditions, where \(u_\ell\) is the lower-precision unit roundoff. On synthetic and real tensors, our methods achieve speedups of a fewfold to several dozenfold over classical HOSVD and STHOSVD, with runtimes competitive with those of other randomized methods. The double-precision versions match the classical reconstruction errors to numerical precision, and the mixed-precision variants yield nearly identical errors.


\end{abstract}

\section{Introduction}\label{sec:introduction}

Higher-order singular value decomposition (HOSVD) is a standard tool for constructing Tucker format of higher-order data~\cite{DLMV2000,KB2009}. Truncated HOSVD computes the dominant left singular subspace of each   unfolding, whereas sequentially truncated HOSVD (STHOSVD) compresses the tensor after processing each mode~\cite{VVM2012}. Both provide well  approximation guarantees, while their use on large tensors can be limited by the cost of computing singular subspaces: Even when each mode dimension is only moderate, the corresponding unfolding may have an enormous other dimension. Our goal is to accelerate classical HOSVD and STHOSVD while preserving their reconstruction accuracy.

Randomized numerical linear algebra offers a natural way to reduce this cost~\cite{HMT2011,MT2020}.  {Minster et al. construct approximate dominant subspaces for HOSVD and STHOSVD by modewise random sampling, while Che et al. employ power iterations to improve subspace extraction}~\cite{MSK2020,CWY2020}.  {Two-stage methods first construct oversampled factors and then truncate the compressed core deterministically; streaming methods recover Tucker approximations from sketches of the factors and core}~\cite{MLB2024,SGLTU2020}.  {The multilinear Nystr\"om method of Bucci and Robol uses two-sided sketches to construct an oblique projection, while RTSMS of Hashemi and Nakatsukasa combines single-mode random compression with least-squares recovery}~\cite{BR2024,HN2023}.

 {Despite various compression and recovery steps, existing methods generally use sketches to approximate the subspace information needed for Tucker recovery.   Standard randomized SVD illustrates this: a sample sized at the target rank plus an   oversampling margin defines the subspace within which the final approximation is computed. Matching the truncated SVD error requires this space to contain an optimal dominant subspace;   directions missed during sampling cannot be recovered by subsequent operations within that space. Small samples can be effective for rapidly decaying spectra, but may leave an accuracy gap for slower decay or high accuracy targets.} Larger samples or power iterations can improve accuracy, but require extra data access, matrix products and sometimes reorthogonalization, which may reduce computational gains~\cite{HMT2011,CWY2020,LLL2026}.

\begin{figure}[t]
 \centering
 \includegraphics[width=\textwidth]{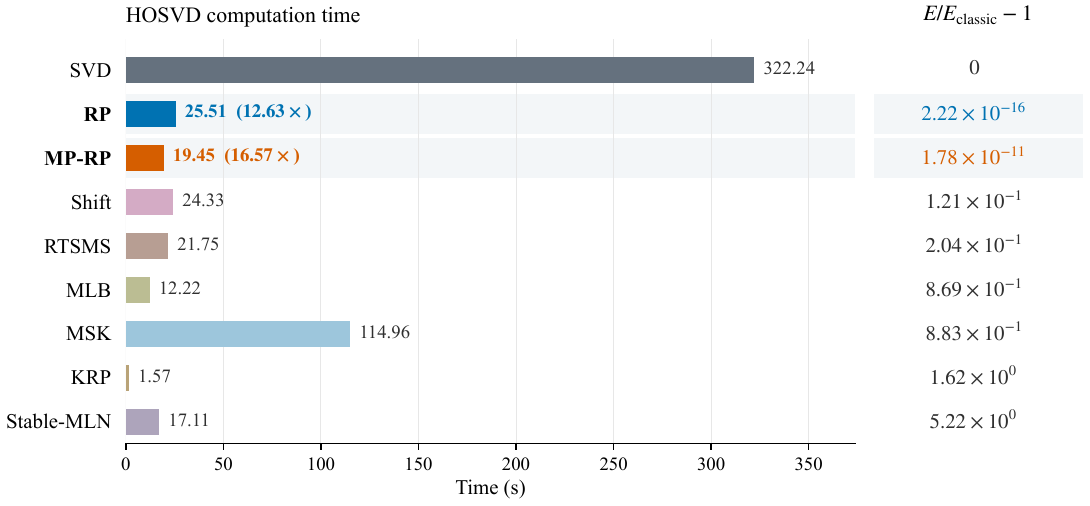}
 \caption{HOSVD tests of a $200^4$ tensor with target ranks $(20,20,20,20)$ and medium spectral decay; see Sect.~\ref{sec:experiment-data} for its construction and Sect.~\ref{sec:experiment-comparison} for the method names. All randomized  methods use sparse sign test matrices. Bars show median times over five runs. The right column reports $E/E_{\mathrm{classic}}-1$, where $E$ is the   Frobenius reconstruction error.  }
 \label{fig:intro-time-accuracy}
 \vspace{-12pt}
\end{figure}

\emph{However},      must random sketches be used primarily to identify approximate dominant subspaces? The shapes of  tensor unfoldings suggest a different possibility. For instance, an unfolding of a cubic  tensor has dimensions of the form \(n\times n^2\); for higher-order tensors, the aspect ratio can be even more extreme. Thus the long dimension may be huge while the mode dimension remains manageable. A sketch whose size scales with the mode dimension can still be dramatically smaller than the original unfolding, even though it is larger than the usual target-rank-scale sample. Furthermore, HOSVD requires only the   dominant singular vectors on the mode side, without needing  its right singular vectors. Together, these observations make it possible to preserve the    entire row space geometry  and recover only the small factor needed by HOSVD.

These observations suggest a different computational role for randomness: instead of   identifying an approximate dominant subspace, random sketches can be used to precondition its recovery from the original unfolding. Motivated by this, we develop randomized preconditioned HOSVD   (RP-HOSVD and RP-STHOSVD). Specifically, we use the singular information of a sketch to construct a preconditioner, numerically truncating near machine precision sketch singular values to avoid unstable inverse scaling. With a suitable oblivious subspace embedding (OSE), the preconditioned matrix has an \(O(1)\) condition number with high probability, independently of the conditioning of the original unfolding. This bounded conditioning enables Gram–Cholesky computations to recover a small factor containing the retained singular information, from which the desired dominant subspace is obtained by target-rank truncation. Preconditioning thus avoids the accuracy risks of directly forming the   Gram matrix of the original unfolding.

The sketching and whitening ingredients have precedents in randomized least-squares preconditioning and randomized Cholesky QR~\cite{RT2008,Blendenpik,LSRN,Balabanov}; our focus is their use for the one-sided singular subspace   required by HOSVD. Connections and differences with existing randomized preconditioning literature will be discussed in Sect. \ref{sec:truncated-rp}.

The computational benefit of the randomized preconditioning differs from that of   small rank-sized sketches. The leading arithmetic complexity may remain of the same order as a direct compact SVD. The opportunity comes from organizing  the large-scale   work as efficient BLAS-3 matrix--matrix operations, while confining the SVD and Cholesky factorizations to matrices governed by the mode dimension. This organization   naturally invites further acceleration via mixed precision: the expensive large matrix products   use    lower precision, while the small factorizations and final recovery remain in higher precision. These variants are termed MP-RP-HOSVD and MP-RP-STHOSVD.

Fig.~\ref{fig:intro-time-accuracy} gives a first impression of the speed and accuracy of the proposed methods on a fourth-order tensor. RP-HOSVD and MP-RP-HOSVD reduce the classical HOSVD time from $322.24$ seconds to $25.51$ and $19.45$ seconds, respectively, corresponding to speedups of $12.63$ and $16.57$. Both methods have runtimes comparable to those of several randomized Tucker methods. RP-HOSVD matches the classical reconstruction error to machine precision, and the relative discrepancy for MP-RP-HOSVD is only $1.78\times10^{-11}$.  Sect.~\ref{sec:numerical-experiments} conducts detailed experiments on synthetic and real tensors.

We analyze   the approximation accuracy and numerical stability of the proposed methods. For RP-HOSVD and RP-STHOSVD, we establish reconstruction bounds that retain the classical HOSVD tail-energy term, with an additional term quadratic in the numerical truncation thresholds and embedding distortion. This   dependence is sharp in general. Moreover,   floating-point analysis shows that, under suitable embedding and threshold conditions, the proposed method   has the same order of normwise backward error as a direct backward-stable SVD, with constants independent of the unfolding's condition number; see Appendix~\ref{app:rp-roundoff}.

For mixed precision, we analyze how errors in the   large matrix products interact with numerical truncation and affect the final projection error. This analysis yields unified bounds for MP-RP-HOSVD and MP-RP-STHOSVD: under suitable embedding and truncation conditions and the rounding model used in Sect.~\ref{sec:mixed-rp}, the squared reconstruction error is bounded by the classical HOSVD tail energy, up to a relative $O(\sqrt{u_\ell})$ correction and an additive $O(u_\ell)$ term. Here $u_\ell$ is the lower-precision unit roundoff.

\paragraph{Code.}
The codes for algorithms, experiment reproduction scripts, and Lean formalizations of the theoretical results are available at \url{https://github.com/yuningyang19/RP_HOSVD}.

The remainder   is organized as follows. Sect.~\ref{sec:preliminaries} introduces notation and subspace embeddings. Sect.~\ref{sec:rp-hosvd} presents the tensor algorithms and reconstruction bounds. Sect.~\ref{sec:rp-matrix} develops and analyzes the randomized preconditioning procedure, and Sect.~\ref{sec:mixed-rp} extends it to mixed precision. Sect.~\ref{sec:numerical-experiments} reports numerical experiments. Sect.~\ref{sec:conclusions} draws conclusions.

\section{  Preliminaries}\label{sec:preliminaries}

\subsection{ Tensor unfoldings and   Tucker approximation}

Let $[D]_r$ denote a best rank-$r$ approximation to a matrix $D$, and   $\tau_{r+1}^2(D)=\sum_{i>r}\sigma_i(D)^2$ the tail energy. Let $\kappa_2(D)=\sigma_{\max}(D)/\sigma_{\min}(D)$ denote the $2$-norm condition number.

Let $\mathcal A\in\R^{n_1\times\cdots\times n_d}$ be an order-$d$ tensor. Denote its mode-$j$ unfolding by $A_j\in\R^{n_j\times M_j}$ with $M_j=\prod_{k\ne j}n_k$. We use a fixed column ordering for each unfolding.
For $W\in\R^{m\times n_j}$,   $\mathcal A\times_j W$  denotes the mode-$j$ product, and its mode-$j$ unfolding is $WA_j$. Products along different modes commute. For these conventions and basic Tucker properties, see~\cite{KB2009}.

Given a target multilinear rank $(r_1,\ldots,r_d)$  with $1\le r_j\le n_j$, one seeks orthonormal factors $U_j\in\R^{n_j\times r_j}$  and the core tensor $\mathcal G=\mathcal A\times_1 U_1\T\cdots\times_d U_d\T$ such that
 \begin{align}\label{eq:tucker-projection}
 \widehat{\mathcal A}
 &=\mathcal G\times_1 U_1\cdots\times_d U_d
 =\mathcal A\times_1(U_1U_1\T)\cdots\times_d(U_dU_d\T)
 \end{align}
has a small reconstruction error.
The main algorithmic choice is   how to obtain the subspace for each mode.

The classical truncated HOSVD~\cite{DLMV2000} takes the leading $r_j$ left singular vectors of each original unfolding $A_j$  and then computes the core. The sequentially truncated HOSVD (STHOSVD) compresses each mode immediately after obtaining its factor, so the next step works with a smaller tensor~\cite{VVM2012}.  For simplicity, the sequential version processes the modes in the order $1,\ldots,d$; any other fixed order can be obtained by permuting the modes. Algorithm~\ref{alg:classical-hosvd} summarizes both procedures.

\begin{algorithm}[htbp]
\caption{Classical truncated HOSVD and STHOSVD}\label{alg:classical-hosvd}
\begin{algorithmic}[1]
\Require   $\mathcal A\in\R^{n_1\times\cdots\times n_d}$, target ranks $(r_1,\ldots,r_d)$, and a choice of HOSVD or STHOSVD.
\Ensure Orthonormal factors $U_j\in\R^{n_j\times r_j}$   and core $\mathcal G$.
\If{HOSVD is selected}
 \For{$j=1,\ldots,d$}
  \State Compute an SVD of $A_j$ and form $U_j$ from its leading $r_j$ left singular vectors.
 \EndFor
 \State $\mathcal G=\mathcal A\times_1 U_1\T\cdots\times_d U_d\T$.
\Else\Comment{STHOSVD, modes $1,\ldots,d$}
 \State $\mathcal A^{(0)}=\mathcal A$.
 \For{$j=1,\ldots,d$}
  \State Take the mode-$j$ unfolding $D_j$ of the current tensor $\mathcal A^{(j-1)}$.
  \State Compute an SVD of $D_j$ and form $U_j$ from its leading $r_j$ left singular vectors.
  \State $\mathcal A^{(j)}=\mathcal A^{(j-1)}\times_j U_j\T$.
 \EndFor
 \State $\mathcal G=\mathcal A^{(d)}$.
\EndIf
\State \Return $\mathcal G,U_1,\ldots,U_d$.
\end{algorithmic}
\end{algorithm}

Both   algorithms satisfy the tail-energy bound~\cite{DLMV2000,VVM2012}

\begin{equation}\label{eq:classical-tucker-bound}
 \norm{\mathcal A-\widehat{\mathcal A}}_F^2
 \le\sum_{j=1}^d\tau_{r_j+1}^2(A_j).
\end{equation}

\subsection{Oblivious subspace embeddings}

Let $\mathscr{U}\subseteq\R^M$ be a fixed subspace with $\dim\mathscr{U}\le k$, and $0<\eta,\delta<1$. Suppose $S\in\R^{s\times M}$ is drawn independently of $\mathscr{U}$ and, for every such subspace, satisfies the following with probability at least $1-\delta$:

\begin{equation}\label{eq:ose-standard}
 (1-\eta)\norm{x}_2^2\le\norm{Sx}_2^2
 \le(1+\eta)\norm{x}_2^2,
 \qquad \text{  for all }x\in\mathscr{U};
\end{equation}
then this distribution is called an $(\eta,\delta)$-oblivious subspace embedding, or shortly $(\eta,\delta)$-OSE, for subspaces of dimension at most $k$. Here $\eta$ is the embedding distortion and $\delta$ bounds the failure probability~\cite{Woodruff2014}.
If $Q$ is an orthonormal basis for $\mathscr{U}$, then~\eqref{eq:ose-standard} is equivalent to
\[
(1-\eta)I\preceq Q^\top S^\top SQ\preceq(1+\eta)I.
\]
 This property controls the condition number of a preconditioned matrix.

Several OSE constructions are available. A standard Gaussian sketch has independent entries $S_{ij}\sim N(0,1/s)$. For any fixed subspace of dimension at most $k$, a sufficient choice for~\eqref{eq:ose-standard} is $s\ge c\eta^{-2}(k+\log(1/\delta))$, where $c$ is an absolute constant~\cite{Woodruff2014}. Thus the embedding dimension need not grow with the ambient dimension $M$. Forming a dense Gaussian sketch directly can still be costly. A sparse sign sketch places only a few random signs in each column, with suitable normalization, reducing the cost of applying $S$. Its embedding dimension and column sparsity must be chosen together according to the distribution and embedding guarantee~\cite{Woodruff2014}. For tensor unfoldings, structured random projections such as Khatri--Rao projections can also use the mode structure to organize the sketch~\cite{SVB2026}. The required sketch dimension depends on the embedding distribution and the prescribed distortion and failure probability.

Throughout this work, the algorithms and analysis   only require $S$ to be an OSE, without prescribing a particular construction.

\section{HOSVDs via Randomized Preconditioning}\label{sec:rp-hosvd}

Our aim is to accelerate the classical truncated HOSVD while preserving its accuracy. The main technique is randomized preconditioning.

\subsection{From randomized range finding to preconditioned recovery}\label{sec:rangefinding-recovery}

We make the difference between subspace capture and preconditioned recovery precise for one mode unfolding \(D\in\R^{n\times M}\), where \(M\gg n\), and a target rank \(r\le n\). A standard rank-sized rangefinder computes an orthonormal basis \(Q\) of a sampled column space $D\Omega$ with $\Omega$ being a random test matrix, typically using a sketch dimension \(s=r+p\), and selects the final rank-\(r\) subspace within \(\range(Q)\)~\cite{HMT2011,MSK2020}. Its associated constrained approximation problem satisfies
\begin{align}\label{eq:rangefinder-residual}
 \min_{\substack{U\T U=I_r\\ \range(U)\subseteq\range(Q)}}
 \norm{D-UU\T D}_F^2
 &=\norm{D-Q[Q\T D]_r}_F^2\nonumber\\
 &=\norm{(I-QQ\T)D}_F^2+\tau_{r+1}^2(Q\T D).
\end{align}
The first term measures information outside the candidate space, which cannot be restored by rank truncation within that space. Hence the result matches a best rank-\(r\) approximation to \(D\) only when the sampled space admits an optimal dominant subspace.

Our construction uses the sketch for preconditioning. Write $H=D\T\in\R^{M\times n}$ and form $Y=SH$, with a sketch dimension of order $n$. We require $S$ to preserve the geometry of the entire $\range(H)$. The singular information of $Y$ provides preconditioned coordinates, in which we revisit $H$ and recover a small factor before the final rank-$r$ truncation. If the sketch preserves rank and all nonzero sketch singular directions are retained, this factor has the same singular values and right singular subspaces as $H$ in exact arithmetic.

If   rounding errors could be ignored, forming $H\T H$ and computing its eigendecomposition might be the fastest way to obtain the   singular subspace required by HOSVD. In practice, however, Gram formation squares the condition number and can cause substantial accuracy loss; see Sect.~\ref{sec:gram-conditioning}.  Randomized preconditioning provides an efficient way to make the matrix well conditioned before Gram formation. Constructing the preconditioner, however, involves inverse scaling by sketch singular values. We therefore truncate very small values to avoid amplifying rounding errors. The effect of discarding these small singular values on the final rank-$r$ approximation error is analyzed in Sect.~\ref{sec:rp-matrix}.

This procedure organizes the large computations as matrix multiplications, while keeping the factorizations small. It also invites further acceleration via mixed precision: the large products use lower precision, while the small factorizations remain in higher precision; see Sect.~\ref{sec:mixed-rp}. We next incorporate these procedures into HOSVD and STHOSVD.

\subsection{ RP-HOSVD and RP-STHOSVD}

We use either the randomized preconditioning (RP) procedure   (Algorithm~\ref{alg:rp}) or its mixed-precision version MP-RP  (Algorithm~\ref{alg:mixed-rp}) to compute the dominant singular subspaces for each mode of the tensors. Both procedures take a thin matrix, target rank $r$, sketch dimension $s$, and numerical   threshold $\epsilon$ as inputs.  $\epsilon$, usually chosen at the machine precision scale,    prevents inverse scaling by     singular values near machine precision from amplifying recovery errors. $\epsilon$  is assumed throughout to be smaller than the  $r$-th singular value of the sketch. MP-RP   takes   two precisions $u_\ell,u_h$.

RP-HOSVD applies the selected procedure to all original unfoldings, then computes the core. RP-STHOSVD proceeds sequentially, compressing each mode as soon as its factor is obtained; see Algorithms~\ref{alg:rp-hosvd} and~\ref{alg:rp-sthosvd}. At step $j$ of Algorithm~\ref{alg:rp-sthosvd}, the unfolding has size $n_j\times M_j^{\mathrm{seq}}$, where $M_j^{\mathrm{seq}}=\prod_{k<j}r_k\prod_{k>j}n_k$. When $M_j^{\mathrm{seq}}$ is close to $n_j$, we may instead compute the SVD of $D_j$ directly.

\begin{algorithm}[htbp]
\caption{RP-HOSVD}\label{alg:rp-hosvd}
\begin{algorithmic}[1]
\Require   $\mathcal A\in \mathbb R^{n_1\times\cdots \times n_d}$, target ranks $(r_j)_{j=1}^d$, sketch dimensions $s_j$, thresholds $\epsilon_j$, and RP or MP-RP; also $u_\ell,u_h$ for MP-RP.

\Ensure Factors $U_j\in\R^{n_j\times r_j}$ and core $\mathcal G$.
\For{$j=1,\ldots,d$}
\Statex Take the original unfolding $A_j$.
\Statex Apply  RP or MP-RP to $A_j\T$ with parameters $(r_j,s_j,\epsilon_j)$ to obtain $U_j$; use the specified $u_\ell,u_h$ if MP-RP is used.
\EndFor
\State $\mathcal G=\mathcal A\times_1 U_1\T\cdots\times_d U_d\T$.
\end{algorithmic}
\end{algorithm}

\begin{algorithm}[htbp]
\caption{RP-STHOSVD}\label{alg:rp-sthosvd}
\begin{algorithmic}[1]
\Require The same inputs as Algorithm~\ref{alg:rp-hosvd}, with mode order $1,\ldots,d$.

 \Ensure Factors $U_j\in\R^{n_j\times r_j}$ and core $\mathcal G$.
\State $\mathcal A^{(0)}=\mathcal A$.
\For{$j=1,\ldots,d$}
 \State Take the mode-$j$ unfolding $D_j$ of $\mathcal A^{(j-1)}$.

 \State  Apply RP or MP-RP to $D_j\T$ with parameters $(r_j,s_j,\epsilon_j)$ to obtain $U_j$; use the specified $u_\ell,u_h$ if MP-RP is used.

 \State $\mathcal A^{(j)}=\mathcal A^{(j-1)}\times_j U_j\T$.
\EndFor
\State $\mathcal G=\mathcal A^{(d)}$.
\end{algorithmic}
\end{algorithm}

\paragraph{ Computational     cost.}

The cost per mode is given by~\eqref{eq:rp-cost};
summing these costs and adding the core contractions gives the total cost. With $s_j=O(n_j)$, RP has the same $O(Mn_j^2)$ leading cost as a direct compact SVD, but gains speed through its computational organization.    RP puts all $O(Mn_j^2)$ work into BLAS-3 matrix multiplications; its factorizations and other BLAS-2 operations involve only $O(n_j)\times n_j$ matrices and cost only $O(n_j^3)$. Efficient BLAS-3 computation and potentially smaller leading constants can therefore give substantial practical gains even with the same leading-order complexity. LSRN~\cite{LSRN} offers a similar example: despite the same $O(Mn^2)$ complexity and possibly more flops, it can outperform a direct SVD for similar reasons.

 Mixed precision   further reduces time at the same complexity by using low precision for   large multiplications and high precision for   small factorizations.

\subsection{Error bounds}

Let $U_j$ be arbitrary orthonormal factors. Denote the original mode-$j$ unfolding by $A_j$ and construct $\widehat{\mathcal A}$ via $U_j$ by~\eqref{eq:tucker-projection}. HOSVD and STHOSVD satisfy~\cite{VVM2012,DLMV2000}

\begin{equation}\label{eq:tensor-local-errors}
 \norm{\mathcal A-\widehat{\mathcal A}}_F^2
 \le\sum_{j=1}^d\norm{(I-U_jU_j\T)A_j}_F^2.
\end{equation}
Below, $D_j$ uniformly denotes the unfolding at step $j$: the original $A_j$ for RP-HOSVD, or the current compressed unfolding for RP-STHOSVD. $\widehat{\mathcal A}$ is constructed by $U_j$ returned by RP-HOSVD or RP-STHOSVD.

\begin{theorem}[RP-HOSVD and RP-STHOSVD] \label{thm:tensor-rp}

Let $0<\eta,\delta<1$. For each step $j$, independently draw an $(\eta,\delta)$-OSE $S_j$ for subspaces of dimension at most $n_j$, and apply Algorithm~\ref{alg:rp} (RP) to $D_j\T$ to obtain $U_j$. Let $\epsilon_j>0$ be small enough to satisfy~\eqref{eq:rp-scale-separation} of Theorem~\ref{thm:rp-refined} at each step. Then, in exact arithmetic, with probability at least $1-d\delta$,

\begin{equation}\label{eq:tensor-rp-refined}
 \norm{\mathcal A-\widehat{\mathcal A}}_F^2
 \le\sum_{j=1}^d\tau_{r_j+1}^2(A_j)
 +\frac{2\eta^2}{1-\eta}\sum_{j=1}^d r_j\epsilon_j^2.
\end{equation}
Especially, if $\epsilon_j$ is smaller than every nonzero singular value of $S_jD_j\T$, then the   $\epsilon$-term in \eqref{eq:tensor-rp-refined}   vanishes, giving the classical tail-energy bound~\eqref{eq:classical-tucker-bound}.

\end{theorem}
\begin{proof}
Assume first that~\eqref{eq:ose-standard} holds at every step with $\mathscr U=\range(D_j\T)$ and $S=S_j$. Theorem~\ref{thm:rp-refined} applied to $D_j\T$ then gives
\begin{equation}\label{eq:tensor-mode-residual}
 \norm{(I-U_jU_j\T)D_j}_F^2
 \le \tau_{r_j+1}^2(D_j)
      +\frac{2\eta^2}{1-\eta}r_j\epsilon_j^2.
\end{equation}
For RP-HOSVD, $D_j=A_j$. \eqref{eq:tensor-mode-residual} together with \eqref{eq:tensor-local-errors}   gives~\eqref{eq:tensor-rp-refined}.

For RP-STHOSVD, as the error decomposition for successive orthogonal projections~\cite[Theorem~5.1]{VVM2012}, one has $\| \mathcal A-\widehat{\mathcal A} \|_F^2=\sum_{j=1}^d\norm{(I-U_jU_j\T)D_j}_F^2$. On the other hand, as in the proof of~\cite[Theorem~6.5]{VVM2012}, one gets $\tau_{r_j+1}(D_j)\le\tau_{r_j+1}(A_j)$. These, together with~\eqref{eq:tensor-mode-residual}, give~\eqref{eq:tensor-rp-refined}.

It remains to show that~\eqref{eq:tensor-rp-refined} holds with probability at least $1-d\delta$. Since $S_j$ is independent of $D_j$, \eqref{eq:ose-standard} holds on $\range(D_j\T)$ with probability at least $1-\delta$. A union bound shows that this condition holds at all steps with probability at least $1-d\delta$.

The last assertion that the $\epsilon$-term in \eqref{eq:tensor-rp-refined} vanishes also can be deduced from Theorem~\ref{thm:rp-refined}.
\end{proof}

Theorem~\ref{thm:tensor-rp} shows that both RP-HOSVD and RP-STHOSVD retain the classical tail-energy term, with only an additional $O(\eta^2\epsilon_j^2)$ term from numerical truncation. With   $\epsilon_j$ being machine precision order, both algorithms nearly recover the classical bound. The example after the proof of Theorem~\ref{thm:rp-refined} shows that the order $\eta^2\epsilon_j^2$ is sharp in general.

Mixed precision further requires bounded entries of $S$ to control rounding errors. Examples include sparse sign sketches~\cite{Woodruff2014}, Khatri--Rao projections of Rademacher matrices~\cite{SVB2026}, and Kronecker-structured SRHTs~\cite{MLB2024}.

\begin{theorem}[Mixed-precision RP-HOSVD and RP-STHOSVD]\label{thm:tensor-mp}

Use    Algorithm \ref{alg:mixed-rp} (mixed-precision RP) and draw independent $(\eta,\delta)$-OSE $S_j$ at each step of RP-HOSVD and RP-STHOSVD,   with the   assumptions of Theorem~\ref{prop:mp-recovery} holding. Suppose that the entries of all $S_j$ are uniformly bounded. Take $\epsilon_j=\Theta\!\left(\sqrt{u_\ell}\norm{D_j}_F\right)$. Then  with probability at least $1-d\delta$, there holds
\begin{equation}\label{eq:tensor-mp-bound}
 \begin{aligned}
 \norm{\mathcal A-\widehat{\mathcal A}}_F^2
 \le{}&\bigl(1+O(\sqrt{u_\ell})\bigr)
             \sum_{j=1}^d\tau_{r_j+1}^2(A_j)+O(u_\ell).
 \end{aligned}
\end{equation}
\end{theorem}
\begin{proof}
By~\eqref{eq:mp-asymptotic-accuracy} in Theorem~\ref{prop:mp-recovery},   at each step there holds
\[
 \norm{(I-U_jU_j\T)D_j}_F^2
 \le\bigl(1+O(\sqrt{u_\ell})\bigr)\tau_{r_j+1}^2(D_j)+O(u_\ell).
\]
The remaining error and probability arguments follow the proof of Theorem~\ref{thm:tensor-rp}.
\end{proof}

 \eqref{eq:tensor-mp-bound} highlights the scale of the errors introduced by low precision: the tail-energy term increases by only a relative $O(\sqrt{u_\ell})$, with an additional $O(u_\ell)$ error. For a more explicit estimate of a single mode, see Theorem~\ref{prop:mp-recovery}.

\section{Computing Left Dominant Subspace via Randomized Preconditioning}\label{sec:rp-matrix}
We   describe the RP procedure used in the previous section. In the sequel, we drop the mode index, write the current unfolding $D_j$ as $D\in\R^{n\times M}$ and set $H=D\T\in\R^{M\times n}$, with $M\gg n$.

\subsection{Gram computation and  conditioning problems} \label{sec:gram-conditioning}

\begin{figure}[!htbp]
 \centering
 \includegraphics[width=\textwidth]{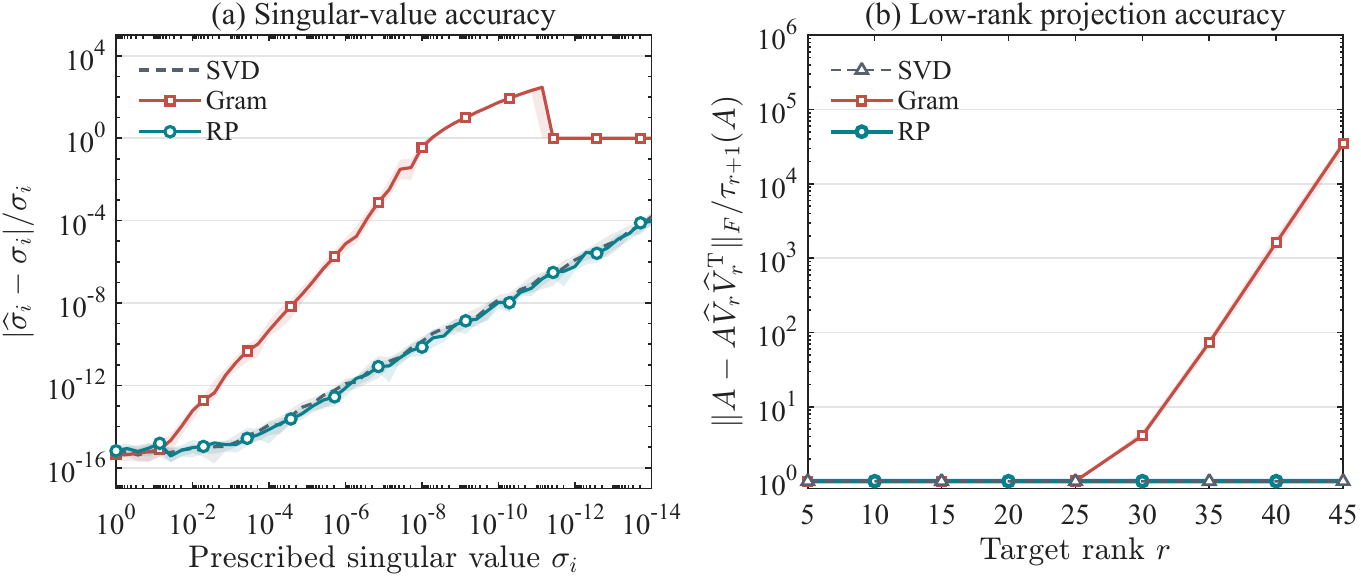}
 \caption{Accuracy for $2000\times50$ matrices with random orthogonal factors and
 prescribed singular values geometrically decreasing from $1$ to $10^{-14}$.
 Left: relative singular-value error
 versus singular values $\sigma_i$.
 Right: relative projection error
versus the target rank $r$, where $\widehat{V}_r\widehat V_r\T$ projects onto the subspace computed by each method.
 SVD, Gram, and RP denote standard compact SVD, Gram eigendecomposition, and Algorithm~\ref{alg:rp}.
 {\color{black}RP uses a sparse sign $S$ with $s=100$ and
 $\epsilon=2\varepsilon_{\mathrm{mach}}\|Y\|_2$, where $\varepsilon_{\mathrm{mach}}$ is MATLAB's machine epsilon.
 }
 Curves and bands show medians and $20$th--$80$th percentiles over $20$   trials.
 Interlaced triangles (SVD) and circles (RP) distinguish their nearly overlapping curves on the right.
}
 \label{fig:gram-rp-accuracy}
\end{figure}

In HOSVD and STHOSVD, one only needs the left dominant subspace of the unfolding (the right counterpart of $H$). A naive approach is to form $H\T H$ and compute its eigendecomposition. This requires only one multiplication involving the long dimension $M$ and a small $n\times n$ factorization. In exact arithmetic, it gives the same dominant subspace as a direct SVD. In floating point, however, as $\kappa_2(H\T H)=\kappa_2(H)^2$, Gram formation may obscure small spectral components that a direct SVD can still resolve. Fig.~\ref{fig:gram-rp-accuracy} illustrates this. On the left, the  singular value errors from Gram formation grow as the singular values decrease, while direct SVD and RP (Algorithm~\ref{alg:rp}) still resolve smaller values. The right shows the  relative projection errors. All three methods are nearly optimal at small ranks, but the Gram residual departs markedly from the optimum at larger ranks; RP stays close to SVD.

In view of the above, our goal is to keep the factorizations small without squaring the condition number of $H$. First assume that $H$ has full column rank. If an invertible $T\in\R^{n\times n}$ makes $X=HT$ well conditioned, we can form its Gram matrix without significant loss of accuracy and compute the Cholesky factor from $R\T R=X\T X$. Setting $B=RT^{-1}\in\mathbb R^{n\times n}$ then gives
\begin{equation}\label{eq:ideal-gram}
 B\T B=T^{-T}R^TRT^{-1}=T^{-\mathsf T}X\T XT^{-1}=H\T H.
\end{equation}
Thus $B$ and $H$ have the same singular values and right singular subspaces.
Therefore,
the remaining task is to cheaply construct $T$.

\subsection{Constructing \texorpdfstring{$T$}{T}}

First assume that $H$ has full column rank. Draw a random test matrix $S\in\R^{s\times M}$ and form the sketch
\begin{equation}\label{eq:sketch}
 Y=SH \in \mathbb R^{s\times n}.
\end{equation}
We consider \(s=cn\), where \(c\ge1\) is a constant, so that \(Y\) is much smaller than \(H\) when \(M\gg n\).
Take $\mathscr U=\range(H)$ in~\eqref{eq:ose-standard} and   assume that $S$ is an $(\eta,\delta)$-OSE, namely, with probability at least $1-\delta$,
\[
 (1-\eta)\norm{Hz}_2^2\le\norm{SHz}_2^2\le (1+\eta)\norm{Hz}_2^2,
 \qquad z\in\R^n.
\]
Henceforth, we fix a realization of $S$ for which the above
embedding bound holds. This embedding property ensures that  $Y=SH$ also has full column rank. Take the compact SVD of $Y$:
\begin{equation}\label{eq:sketch-svd}
 Y=U\Sigma V\T,~~~~ U\T U = I_n, V\T V = I_n,
\end{equation}
where  $\Sigma$ has positive diagonal entries in decreasing order. Set
\[
 T=V\Sigma^{-1},\qquad X=HV\Sigma^{-1}.
\]
To estimate the condition number of $X$, note that
\[
 SX=SHV\Sigma^{-1}=YV\Sigma^{-1}=U;
\]
Hence $\norm{SXy}_2=\norm{y}_2$ for every $y\in\R^n$. Since $Xy\in\range(H)$, the embedding property~\eqref{eq:ose-standard} gives $(1-\eta)\norm{Xy}_2^2\le\norm{SXy}_2^2=\norm{y}_2^2\le(1+\eta)\norm{Xy}_2^2$. Thus when the distortion is small,   $ X $ is well conditioned.   This use of sketch-based coordinates is standard  \cite{RT2008,LSRN}.

Now, form $R\T R=X\T X$ and $B=R\Sigma V\T$. By~\eqref{eq:ideal-gram},   the leading $r$ right singular vectors of $B$ give a best rank-$r$ projection of $H$ (the left counterpart of $D$).   Randomization serves only to construct the preconditioner, a different role from that in the randomized HOSVD methods  \cite{MSK2020,CWY2020,SGLTU2020,MLB2024,LLL2026,SVB2026}.

\subsection{Dealing with small singular values} \label{sec:truncated-rp}

Forming $T=V\Sigma^{-1}$ can amplify errors in the presence of too small singular
values, and is not applicable when $H$ is rank
deficient. We therefore retain and invert only singular values above
a prescribed threshold.

We allow rank-deficient $H$.
Let $Y$ have   SVD~\eqref{eq:sketch-svd}, with
singular values in decreasing order, and choose a threshold
$\epsilon\ge0$. In practice, one may take $\epsilon=cu\norm{Y}_2$,
where $u$ is the unit roundoff and $c$ is a constant,
to avoid inverting singular values near the rounding-error
scale while keeping truncation loss small. The analysis below
applies to any  $\epsilon\geq 0$. Let
\begin{equation}\label{eq:retained}
 t=\#\{i:\sigma_i(Y)>\epsilon\},~~
 \Sigma_\epsilon=\diag(\sigma_1(Y),\ldots,\sigma_t(Y)),~~
 V_\epsilon=V(:,1:t),
\end{equation}
and let $U_\epsilon$ denote the matrix of the corresponding left singular vectors.
Define
\begin{equation}\label{eq:whitening}
 T_\epsilon=V_\epsilon\Sigma_\epsilon^{-1}\in\R^{n\times t},\qquad X_\epsilon=HT_\epsilon\in\R^{M\times t}.
\end{equation}
Since $YV_\epsilon=U\Sigma V^TV_\epsilon=U_\epsilon\Sigma_\epsilon$, we have
\begin{equation}\label{eq:sketched-orthogonality}
 SX_\epsilon=YV_{\epsilon}\Sigma_\epsilon^{-1}=\ U_\epsilon,
\end{equation}
We have the following condition number bound.
\begin{lemma}\label{lem:conditioning}
Let $t\ge1$ and let $X_\epsilon$ be defined by~\eqref{eq:whitening}. If~\eqref{eq:ose-standard} holds with $\mathscr U=\range(H)$, then $X_\epsilon$ has full column rank and
\begin{equation}\label{eq:conditioning}
 \frac{1}{1+\eta}I_t\preceq X_\epsilon\T X_\epsilon\preceq \frac{1}{1-\eta}I_t,
 \qquad \kappa_2(X_\epsilon)\le\sqrt{\frac{1+\eta}{1-\eta}}.
\end{equation}
\end{lemma}

\begin{proof}
By~\eqref{eq:sketched-orthogonality}, $\norm{SX_\epsilon y}_2=\norm{y}_2$ for every $y\in\R^t$. Since $X_\epsilon y\in\range(H)$, applying~\eqref{eq:ose-standard} gives
\(
 (1-\eta)\norm{X_\epsilon y}_2^2
 \le\norm{SX_\epsilon y}_2^2
 =\norm{y}_2^2
 \le (1+\eta)\norm{X_\epsilon y}_2^2
\).
Rearranging the two inequalities gives
\[
 \frac{1}{1+\eta}\norm{y}_2^2\le\norm{X_\epsilon y}_2^2
 \le \frac{1}{1-\eta}\norm{y}_2^2.
\]
This proves~\eqref{eq:conditioning} and the full column rank of $X_\epsilon$.
\end{proof}

For moderate distortion $\eta$, the lemma gives
$\kappa_2(X_\epsilon)=O(1)$, independently of $\kappa_2(H)$.
This enables Gram--Cholesky and subsequent   computations:
\begin{equation}\label{eq:truncated-B}
 R_\epsilon\T R_\epsilon=X_\epsilon\T X_\epsilon,\qquad
 B_\epsilon=R_\epsilon\Sigma_\epsilon V_\epsilon\T
 \in\R^{t\times n}.
\end{equation}
If $t=n$, then $B_\epsilon=B$, whose SVD gives the dominant
right singular subspace of $H$. With truncation, $B_\epsilon$
instead recovers the Gram   of $HP_\epsilon$, where
$P_\epsilon=V_\epsilon V_\epsilon\T$. Indeed, \eqref{eq:whitening} gives
$T_\epsilon\Sigma_\epsilon V_\epsilon\T=P_\epsilon$ and
$X_\epsilon=HT_\epsilon$, so~\eqref{eq:truncated-B} yields
\begin{equation}\label{eq:truncated-B-gram}
 B_\epsilon\T B_\epsilon
 =V_\epsilon\Sigma_\epsilon X_\epsilon\T X_\epsilon
   \Sigma_\epsilon V_\epsilon\T = (V_\epsilon\Sigma_\epsilon T_\epsilon\T)H^T H (T_\epsilon\Sigma_\epsilon V_\epsilon\T)
 =(HP_\epsilon)\T(HP_\epsilon).
\end{equation}
Thus $B_\epsilon$ preserves the   right singular information of $HP_\epsilon$.   With a machine-precision threshold $\epsilon$, $HP_\epsilon\approx H$ and $B_\epsilon\T B_\epsilon\approx H\T H$. Thus we  use the leading right singular vectors of $B_\epsilon$
to approximate those of $H$.
Algorithm~\ref{alg:rp} summarizes the procedure, and
Proposition~\ref{thm:rp} bounds the projection error.

\begin{algorithm}[htbp]
\caption{Randomized preconditioned SVD   (RP)}\label{alg:rp}
\begin{algorithmic}[1]
\Require   $H\in\R^{M\times n}$ with $M\gg n$, target rank $ r $, sketch dimension $s$, and threshold $\epsilon\ge0$.
\Ensure  Orthonormal factor $\widehat V_r\in\R^{n\times r}$.
\State  Draw a random test matrix $S\in\R^{s\times M}$ and compute $Y=SH$.
\State
 Compute the SVD $Y=U\Sigma V\T$ and form $t$,
$\Sigma_\epsilon$, and $V_\epsilon$ as in~\eqref{eq:retained},
retaining only singular values greater than $\epsilon$.
\State   Compute $X_\epsilon=HV_\epsilon\Sigma_\epsilon^{-1}\in\R^{M\times t}$.
\State   Form $G_\epsilon=X_\epsilon\T X_\epsilon$ and compute its upper triangular Cholesky factor $R_\epsilon$.
\State   Form $B_\epsilon=R_\epsilon\Sigma_\epsilon V_\epsilon\T$
and take $\widehat V_r$ as its leading $r$ right singular vectors.
\State \Return $\widehat V_r$.
\end{algorithmic}
\end{algorithm}

\paragraph{Stability.} Under the embedding assumptions of Proposition~\ref{prop:rp-backward-stability}, with $\epsilon=cu\norm{Y}_2$ and sufficiently large fixed $c$, the recovered singular values and right singular vectors are those of $H+\Delta H$, where $\norm{\Delta H}_F\le Cu\norm{H}_F$ and $C$ is independent of $\kappa_2(H)$.

\paragraph{Computational cost.}
For dense unstructured $H$, let $T_{\mathrm{sk}}(H,S)$ be the cost of forming $Y=SH$. For instance, a sparse $S$ with at most $\zeta$ nonzeros per column costs $O(\zeta Mn)$. The SVD of $Y$ costs $O(sn^2)$. Forming $X_\epsilon$, its Gram matrix, and the Cholesky factor costs $O(Mnt)$, $O(Mt^2)$, and $O(t^3)$, respectively; forming $B_\epsilon$ and computing its SVD costs $O(nt^2)$. The total is
\begingroup
\abovedisplayskip=5pt
\belowdisplayskip=5pt
\abovedisplayshortskip=3pt
\belowdisplayshortskip=5pt
\begin{equation}\label{eq:rp-cost}
 T_{\mathrm{sk}}(H,S)+O(sn^2+Mnt+Mt^2+nt^2+t^3).
\end{equation}
\endgroup
For $t=n$, RP has the same leading cost $O(Mn^2)$ as a compact
SVD of $H$. Nevertheless, all the $O(Mn^2)$ work consists of efficient BLAS-3
multiplications forming  $Y$, $X_\epsilon$, and $G_\epsilon$,
with only $O(n)\times O(n)$ matrices to factorize. When $M\gg n$,  this makes RP
substantially faster in practice.

\paragraph{Relation to randomized preconditioning literature.}
Randomized preconditioning was used in the  
least-squares algorithm and Blendenpik~\cite{RT2008,Blendenpik}.
Our construction of $T_\epsilon=V_\epsilon\Sigma_\epsilon^{-1}$
from the sketch SVD is similar to that in LSRN~\cite{LSRN},
but the subsequent steps and purpose differ: LSRN iteratively
solves a preconditioned least-squares problem, whereas we
form the Gram   and its Cholesky factor, recover
$B_\epsilon$, and compute its SVD to obtain the factor required by HOSVD.
Our Gram--Cholesky computation is also related to randomized
Cholesky QR~\cite{Balabanov,FGL2021}, but our goal is singular subspace
extraction instead of a QR factorization: we 
use $R_\epsilon\Sigma_\epsilon V_\epsilon\T$ directly, without
explicitly forming the orthogonal factor
$X_\epsilon R_\epsilon^{-1}$. The treatment of rank deficiency also differs: Algorithm~7
in~\cite{Balabanov} selects columns of the original matrix
using strong RRQR, whereas we retain the   singular
vectors of the sketch associated with singular values above
$\epsilon$. 
Balabanov, Melnichenko, and Mahoney~\cite{BMM2026}
use a sketch-QR preconditioner to compute a Q-less QR
factorization of a full column rank operator;
our RP instead uses a sketch-SVD preconditioner, along with numerical
truncation to avoid unstable inverse scaling in rank-deficient
or numerically rank-deficient cases, which were not considered
in~\cite{BMM2026}. Beyond these technical differences, the methods above differ
in purpose from our aim of recovering an exact dominant
right singular subspace of a tall  matrix.
Palitta and Simoncini~\cite{PS2026} also use sketch-based
preconditioning for SVD, but use the sketch's  
    singular vectors as approximations to those of
the original matrix.
  Finally, 
we further establish a refined bound on the excess rank-$r$
projection error caused by numerical truncation, which is
not provided in the works above.

\subsection{Projection error}
We now bound the   error from numerical truncation
and the final rank-$r$ truncation relative to a direct SVD.
\begin{proposition} \label{thm:rp}
Write $P_\epsilon=V_\epsilon V_\epsilon\T$ and $\widehat P_r=\widehat V_r\widehat V_r\T$, where $\widehat V_r$ is returned by Algorithm~\ref{alg:rp}. Suppose that~\eqref{eq:ose-standard} holds on $\mathscr U=\range(H)$. Then the projection error has the exact decomposition
\begin{equation}\label{eq:exact-error}
 \norm{H-H\widehat P_r}_F^2
 =\norm{H(I-P_\epsilon)}_F^2
  +\norm{HP_\epsilon-[HP_\epsilon]_r}_F^2.
\end{equation}
Moreover,
 \begin{equation}\label{eq:tail-bound}
 \norm{H-H\widehat P_r}_F^2
 \le \tau_{r+1}^2(H)
     +\frac{n-t}{1-\eta}\epsilon^2.
\end{equation}
\end{proposition}
\eqref{eq:tail-bound} quantifies the cost of numerical truncation: beyond the optimal   $\tau_{r+1}^2(H)$, it adds only $O(\epsilon^2)$. A machine-precision threshold therefore keeps the error close to  $\tau_{r+1}^2(H)$. If all $\sigma_i(Y)>\epsilon$, then $t=n$ and no truncation occurs: $\widehat P_r$ is a best rank-$r$ projector for $H$, and the right-hand side of \eqref{eq:tail-bound} leaves only $\tau_{r+1}^2(H)$.

\begin{proof}  We first prove two identities $P_\epsilon\widehat P_r=\widehat P_r$ and $H\widehat P_r=[HP_\epsilon]_r$.

 By the definition of $t$ in \eqref{eq:retained}, $\Sigma_\epsilon$ is invertible, and $X_\epsilon$ has full column rank by Lemma~\ref{lem:conditioning}. Hence $\range(B_\epsilon\T)=\range(V_\epsilon)$. Since   $\widehat V_r$ consists of  leading right singular vectors of $B_{\epsilon}$,   $\range(\widehat V_r)\subseteq\range(V_\epsilon)$, giving $P_\epsilon\widehat P_r=\widehat P_r$ and so

\begin{align}\label{eq:AhatP_r}
 H\widehat P_r=HP_\epsilon\widehat P_r=[HP_\epsilon]_r.
\end{align}
where the last equality follows because $\widehat V_r$   consists of dominant right singular vectors of $HP_\epsilon$ (by~\eqref{eq:truncated-B-gram}, $HP_{\epsilon}$ shares the same nonzero singular values and right singular subspaces with $B_\epsilon$). If the best rank-$r$ approximation is not unique, we take
$[HP_\epsilon]_r=H\widehat P_r$.

 Now split   $H(I-\widehat P_r)=H(I-P_\epsilon)+H(P_\epsilon-\widehat P_r)$. Since $P_\epsilon\widehat P_r=\widehat P_r$ implies $(P_\epsilon-\widehat P_r)(I-P_\epsilon)=0$, we have
\begin{align*}
  \left\langle
               H(P_\epsilon-\widehat P_r), H(I-P_\epsilon)\right\rangle  &= \tr\!\left((P_\epsilon-\widehat P_r)H\T H (I-P_\epsilon)
                      \right)\\
 & =\tr\!\left((I-P_\epsilon)(P_\epsilon-\widehat P_r)H\T H\right)=0.
\end{align*}
This together with $H\widehat P_r=[HP_\epsilon]_r$  proves~\eqref{eq:exact-error}.
The second term of \eqref{eq:exact-error} can be bounded by $\tau_{r+1}^2(H)$:
\begin{align*}
 \norm{HP_\epsilon-[HP_\epsilon]_r}_F
 &\le\norm{HP_\epsilon-[H]_rP_\epsilon}_F\\
 &=\norm{(H-[H]_r)P_\epsilon}_F
 \le\norm{H-[H]_r}_F = \tau_{r+1}(H).
\end{align*}
For the first term, by~\eqref{eq:ose-standard},
\[
 (1-\eta)\norm{H(I-P_\epsilon)}_F^2
 \le\norm{SH(I-P_\epsilon)}_F^2
 =\norm{Y(I-P_\epsilon)}_F^2.
\]
As $P_{\epsilon}=V_{\epsilon}V_{\epsilon}\T$, and by the definition of $V_{\epsilon}$ in \eqref{eq:retained},
the right side is the discarded sketch tail energy $\tau_{t+1}^2(Y)$. The definition of $t$ gives $\sum_{j>t}\sigma_j(Y)^2\le(n-t)\epsilon^2$. Substituting both estimates into~\eqref{eq:exact-error} proves~\eqref{eq:tail-bound}.
\end{proof}

We can sharpen~\eqref{eq:tail-bound} further.
\begin{theorem}\label{thm:rp-refined}
Under the assumptions of Proposition~\ref{thm:rp}, suppose
\begin{equation}\label{eq:rp-scale-separation}
 \epsilon^2\le\frac{1-\eta}{2}\sigma_r^2(B_\epsilon),
\end{equation}
 where $B_\epsilon$ was defined in~\eqref{eq:truncated-B}. Then
\begin{equation}\label{eq:rp-refined-bound}
 \norm{H-H\widehat P_r}_F^2
 \le\tau_{r+1}^2(H)
 +\frac{2r\eta^2}{1-\eta}\epsilon^2.
\end{equation}
In particular, if   $\epsilon<\sigma_{\min}(Y) $, then the second term can be removed.
\end{theorem}

Condition~\eqref{eq:rp-scale-separation} separates the numerical threshold from the target singular value of $B_{\epsilon}$. With a machine-precision threshold, and when the spectrum decays not so fast, it is usually easy to satisfy.

The refined bound~\eqref{eq:rp-refined-bound} reveals a further benefit of a better embedding: the numerical error decreases quadratically as $\eta$ tends to zero, a feature not captured by~\eqref{eq:tail-bound}. The example after the proof shows that the order $\eta^2\epsilon^2$ is sharp in general.

To prove the theorem, we first give a technical lemma.

\begin{lemma}\label{lem:retained-subspace-error}
Let $H\in\R^{m\times n}$, let $V\in\R^{n\times t}$ be orthonormal, and suppose that $HV$ has full column rank and $1\le r\le t$. Complete $V$ to an orthogonal matrix $[V,V_\perp]$, and let $Q$ be an orthonormal basis for $\range(HV)$. If
\[
 \norm{HV_\perp}_2^2\le\frac12\sigma_r^2(HV),
\]then
\begin{equation}\label{eq:retained-subspace-error}
 \norm{H-[HV]_rV\T}_F^2
 \le\tau_{r+1}^2(H)+2r\norm{Q\T HV_\perp}_2^2.
\end{equation}
\end{lemma}

\begin{proof}
For $t=n$, $V$ is orthogonal and the result follows directly from best rank-$r$ approximation. Assume $t<n$. We first   compare the leading $r$   singular values of $H$ and $HV$. Complete $Q$ to an orthogonal basis $[Q,Q_\perp]$. In this basis,
\begin{align}\label{eq:block_QAAQ}
 [Q,Q_\perp]\T HH\T[Q,Q_\perp]
 =\begin{bmatrix}K+ZZ\T&ZW\T\\WZ\T&WW\T\end{bmatrix},
\end{align}
where
\begin{equation}\label{eq:retained-block-definitions}
 K=Q\T HVV\T H\T Q,\qquad
 Z=Q\T HV_\perp,\qquad W=Q_\perp\T HV_\perp.
\end{equation}
The definition of $W$ together with the hypothesis gives $\norm{W}_2^2\le\norm{HV_\perp}_2^2\le\sigma_r^2(HV)/2$. Hence
\[
 W\T\bigl(\sigma_r^2(HV)I-WW\T\bigr)^{-1}W\preceq I.
\]
Since $\sigma_r^2(HV)I-WW\T\succ0$, the Schur complement
criterion gives
\begin{align}\label{eq:schur}
 \begin{bmatrix}
 ZZ\T&-ZW\T\\-WZ\T&\sigma_r^2(HV)I-WW\T
 \end{bmatrix}\succeq0.
\end{align}
Combining~\eqref{eq:schur} with~\eqref{eq:block_QAAQ} gives

\begin{align}\label{eq:block_QAAQ-upper-bound}
 \begin{bmatrix}K+ZZ\T&ZW\T\\WZ\T&WW\T\end{bmatrix}
 \preceq
 \begin{bmatrix}K+2ZZ\T&0\\0&\sigma_r^2(HV)I\end{bmatrix}.
\end{align}
Since $HV$ has full column rank, this together with the definition of $Q$ and \eqref{eq:retained-block-definitions} gives $K\succ0$ and $\lambda_i(K)=\sigma_i^2(HV)$, $1\le i\le t$. Thus $\lambda_r(K+2ZZ\T)\ge\lambda_r(K)=\sigma_r^2(HV)$, so the leading $r$ eigenvalues on the right of \eqref{eq:block_QAAQ-upper-bound} come from the upper-left block. By~\eqref{eq:block_QAAQ-upper-bound}, \eqref{eq:block_QAAQ},   and $2ZZ\T\preceq2\norm{Z}_2^2I$,
\begin{align}\label{eq:retained-singular-value-bound}
 \sum_{i=1}^r\sigma_i^2(H)
 \le\sum_{i=1}^r\lambda_i(K+2ZZ\T)
 \le\sum_{i=1}^r\sigma_i^2(HV)+2r\norm{Z}_2^2.
\end{align}

Finally, write $H-[HV]_rV\T=(HV-[HV]_r)V\T+HV_\perp V_\perp\T$. Orthogonality gives
\begin{align*}
 \norm{H-[HV]_rV\T}_F^2
 &=\norm{HV-[HV]_r}_F^2+\norm{HV_\perp}_F^2\\
 &=\norm{HV}_F^2-\sum_{i=1}^r\sigma_i^2(HV)+\norm{HV_\perp}_F^2\\
 &=\norm{H}_F^2-\sum_{i=1}^r\sigma_i^2(HV)\\
 &\le\norm{H}_F^2-\sum_{i=1}^r\sigma_i^2(H)+2r\norm{Q\T HV_\perp}_2^2\\
 &=\tau_{r+1}^2(H)+2r\norm{Q\T HV_\perp}_2^2.
\end{align*}
The inequality uses~\eqref{eq:retained-singular-value-bound} and $Z=Q\T HV_\perp$.
\end{proof}

\begin{proof}[Proof of Theorem~\ref{thm:rp-refined}]
Complete $V_\epsilon$ to an orthogonal matrix
$[V_\epsilon,V_\perp]$ using the discarded right singular
vectors of $Y$. By~\eqref{eq:ose-standard} and the truncation
rule,
\begin{equation}\label{eq:discarded-spectral-bound}
 \norm{HV_\perp}_2^2
 \le\frac{1}{1-\eta}\norm{SHV_\perp}_2^2 = \frac{1}{1-\eta}\norm{YV_\perp}_2^2
 \le\frac{\epsilon^2}{1-\eta}.
\end{equation}
The last inequality uses the threshold $\epsilon$. If $HV_\perp=0$, then $HP_\epsilon=H$ and the result is immediate. Assume $HV_\perp\ne0$.

We first bound the coupling of the retained and discarded parts in the original space. We first have $(SHV_\epsilon)\T SHV_\perp=V_\epsilon\T Y\T YV_\perp=0$. Thus unit vectors $x\in\range(HV_\epsilon)$ and $y\in\range(HV_\perp)$ satisfy $(Sx)\T Sy=0$ and $\norm{S(x+y)}_2^2=\norm{S(x-y)}_2^2$. Replacing $y$ by $-y$ if necessary, assume $x\T y\ge0$. Since $x+y,x-y\in\mathscr U=\range(H)$, applying \eqref{eq:ose-standard} yields
\[
 (1-\eta)\norm{x+y}^2
 \le\norm{S(x+y)}_2^2
 =\norm{S(x-y)}_2^2
 \le (1+\eta)\norm{x-y}^2.
\]
Rearranging gives $|x\T y|\le\eta$. Let $Q_1$ be an orthonormal basis for $\range(HV_\epsilon)$. For any $z\in\R^{n-t}$, orthonormality gives
\[
 \norm{Q_1\T HV_\perp z}_2
 =\sup_{\substack{x\in\range(HV_\epsilon)\\\norm{x}_2=1}}
 |x\T HV_\perp z|.
\]
If $HV_\perp z\ne0$, take $y=HV_\perp z/\norm{HV_\perp z}_2$ in the angle estimate above to obtain
\[
 \norm{Q_1\T HV_\perp z}_2
 \le\eta\norm{HV_\perp z}_2.
\]
This also holds when $HV_\perp z=0$. Taking the supremum over $\norm{z}_2=1$ and using~\eqref{eq:discarded-spectral-bound} gives
\begin{equation}\label{eq:retained-discarded-coupling}
 \norm{Q_1\T HV_\perp}_2^2
 \le\eta^2\norm{HV_\perp}_2^2
 \le\frac{\eta^2}{1-\eta}\epsilon^2.
\end{equation}

We apply Lemma~\ref{lem:retained-subspace-error} with $V=V_\epsilon$ and $Q=Q_1$ and first   check the assumptions of the lemma.
Lemma~\ref{lem:conditioning} gives full column rank of $X_\epsilon$ and hence of $HV_\epsilon=X_\epsilon\Sigma_\epsilon$. Also, \eqref{eq:discarded-spectral-bound}, \eqref{eq:rp-scale-separation}, and $\sigma_r(HV_\epsilon)=\sigma_r(B_\epsilon)$ give
\[
 \norm{HV_\perp}_2^2\le\frac{\epsilon^2}{1-\eta}
 \le\frac12\sigma_r^2(HV_\epsilon).
\]
The lemma therefore applies. By~\eqref{eq:AhatP_r}, $H\widehat P_r=[HP_\epsilon]_r=[HV_\epsilon]_rV_\epsilon\T$, with singular directions chosen as in the algorithm.   Lemma~\ref{lem:retained-subspace-error} and~\eqref{eq:retained-discarded-coupling} yield
\[
 \norm{H-H\widehat P_r}_F^2
 \le\tau_{r+1}^2(H)+2r\norm{Q_1\T HV_\perp}_2^2
 \le\tau_{r+1}^2(H)+\frac{2r\eta^2}{1-\eta}\epsilon^2.
\]
Finally, if all $\sigma_i(Y)>\epsilon$, then $t=n$ and $\widehat P_r$ is a best rank-$r$ projector for $H$, leaving only $\tau_{r+1}^2(H)$ on the right.
\end{proof}

\paragraph{Sharpness of $\Theta(\eta^2\epsilon^2)$.}
\begingroup
\small
 Take $0<\eta<1$ and $h>0$ with $h^2\le1/[2(1+\eta)]$, and set

\[
 H=\begin{bmatrix}1&\eta h\\0&h\sqrt{1-\eta^2}\end{bmatrix},\qquad
 S=\begin{bmatrix}\sqrt{1-\eta^2}&-\eta\\0&1\end{bmatrix},\qquad
 \epsilon=h\sqrt{1-\eta^2}.
\]
Then $Y=SH=\diag(\sqrt{1-\eta^2},\epsilon)$. For $r=1$, the truncation rule gives $P_\epsilon=\widehat P_1=\diag(1,0)$ and $\sigma_1(B_\epsilon)=1$. Since $S\T S$ has eigenvalues $1\pm\eta$ and $\epsilon^2\le(1-\eta)/2$, both~\eqref{eq:ose-standard} and~\eqref{eq:rp-scale-separation} hold. Direct calculation gives
\[
 \norm{H-H\widehat P_1}_F^2-\tau_2^2(H)
 =\frac{2\eta^2h^2}{\sqrt{(1-h^2)^2+4\eta^2h^2}+1-h^2}.
\]
For $0<\eta\le\eta_0<1$, the right side is $\Theta(\eta^2\epsilon^2)$ with uniform constants.
\endgroup

\section{Mixed-Precision Enhancement}\label{sec:mixed-rp}

The cost  of Algorithm~\ref{alg:rp} makes mixed precision particularly attractive: Most work lies in the large BLAS-3 multiplications that form $Y,X_{\epsilon}, G_{\epsilon}$. These are natural candidates for low precision. The relatively small factorizations, on the other hand, determine the preconditioning directions and the recovered dominant subspace, so we keep them in high precision.   Write the two unit roundoffs as $u_\ell$ and $u_h$, with $u_h\ll u_\ell$.

\subsection{The mixed-precision algorithm}

Tildes distinguish mixed-precision matrices from the exact-arithmetic matrices of the preceding section. Detailed steps are given in Algorithm \ref{alg:mixed-rp}.

\begin{algorithm}[htbp]
\caption{Mixed-precision randomized preconditioned SVD (MP-RP)}\label{alg:mixed-rp}
\begin{algorithmic}[1]
\Require   $H\in\R^{M\times n}$, target rank $1\le r\le n$, sketch dimension $s$, threshold $\epsilon\ge0$, and precisions $u_\ell,u_h$.
\Ensure  $\widetilde V_r\in\R^{n\times r}$ stored in high precision.
\State   Draw a random test matrix $S\in\R^{s\times M}$.
\State Compute $\widetilde Y=\operatorname{fl}_\ell(SH)$. \Comment{$u_\ell$: low precision}
\State Compute an SVD of $\widetilde Y$. \Comment{$u_h$: high precision}
\State Retain singular values greater than $\epsilon$, obtaining $t$, $\widetilde\Sigma_\epsilon$, and $\widetilde V_\epsilon$.
\State Form $\widetilde T_\epsilon=\widetilde V_\epsilon\widetilde\Sigma_\epsilon^{-1}$. \Comment{$u_h$: high precision}
\State Compute $\widetilde X_\epsilon=\operatorname{fl}_\ell(H\widetilde T_\epsilon)$. \Comment{$u_\ell$: low precision}
\State Compute $\widetilde G_\epsilon=\operatorname{fl}_\ell(\widetilde X_\epsilon\T\widetilde X_\epsilon)$. \Comment{$u_\ell$: low precision}
\State Use one triangle of $\widetilde G_\epsilon$ to form a symmetric matrix and convert it to high precision.
\State Compute the Cholesky factor $\widetilde R_\epsilon$. \Comment{$u_h$: high precision}
\State Form $\widetilde B_\epsilon=\widetilde R_\epsilon\widetilde\Sigma_\epsilon\widetilde V_\epsilon\T$. \Comment{$u_h$: high precision}
\State Compute and return the leading $r$ right singular vectors $\widetilde V_r$ of $\widetilde B_\epsilon$. \Comment{$u_h$: high precision}
\end{algorithmic}
\end{algorithm}

\paragraph{Relation to existing mixed-precision methods.}
Mixed precision has been studied for randomized range finding and power iterations~\cite{CHP2022}, random projections in SVD/HOSVD~\cite{OY2023}, and iterative refinement of matrix and tensor approximations~\cite{BKMR2025}. Carson and Dau\v{z}ickait\.{e}~\cite{CD2025} study low-precision sketches for preconditioning and iterative refinement. More closely related, Garrison, Chen, and Ipsen~\cite{GCI2026} combine mixed precision, randomized preconditioning, and Cholesky--QR. Their goal is an explicit $Q$ with double-precision orthogonality, so they still apply the preconditioner and perform the subsequent Cholesky--QR in double precision.

Extracting only a   one-sided singular subspace
leaves more room for low precision: we need not recover the
tall matrix's full orthogonal factor.
Algorithm~\ref{alg:mixed-rp} uses low precision for all three
large multiplications, keeping only the relatively small factorizations
and recovery in high precision.
This precision assignment and recovery goal distinguish it
from the methods above. Our analysis quantifies how the three
multiplication errors and numerical truncation jointly affect
the final rank-$r$ projection residual.

\subsection{Projection error}

To make the effect of the three low-precision multiplications
more explicit, we include rounding errors in sketching, preconditioning, Gram formation, and input conversions, but treat the high-precision small-matrix operations as exact.   A related simplification appears in Carson and Dau\v{z}ickait\.{e}~\cite[\S2.1]{CD2024Nystrom}, who assume an exact final SVD in their mixed-precision Nystr\"om analysis while accounting for the other rounding errors.

Write
\begin{equation}\label{eq:mp-errors}
 \widetilde Y=SH+E_s,\qquad
 X_\circ=H\widetilde T_\epsilon,\qquad
 \widetilde X_\epsilon=X_\circ+E_x,\qquad
 \widetilde G_\epsilon=\widetilde X_\epsilon\T\widetilde X_\epsilon+E_g.
\end{equation}
Here $E_s,E_x,E_g$ denote rounding errors in the three products.
$E_s$ affects direction selection; $E_x,E_g$ affect recovery
in the retained subspace.

Use the standard floating-point model and multiplication bound~\cite{Higham2002}, with no overflow or underflow, and write $\gamma_j^{(\ell)}=ju_\ell/(1-ju_\ell)$ for $ju_\ell<1$. We use $\operatorname{fl}_\ell(BC)$ for multiplication after converting both inputs to low precision. Including the conversions, a product with inner dimension $k$ satisfies $|\operatorname{fl}_\ell(BC)-BC|\le\gamma_{k+2}^{(\ell)}|B||C|$, with entrywise absolute values. The three multiplications therefore satisfy

\begin{equation}\label{eq:mp-rounding-bounds}
 \begin{aligned}
 \norm{E_s}_2&\le d_s:=\gamma_{M+2}^{(\ell)}\norm{|S||H|}_2,\\
 \norm{E_x}_2&\le d_x:=\gamma_{n+2}^{(\ell)}\norm{|H||\widetilde T_\epsilon|}_2,\\
 \norm{E_g}_2&\le d_g:=\gamma_{M+2}^{(\ell)}
 \norm{|\widetilde X_\epsilon|\T|\widetilde X_\epsilon|}_2.
 \end{aligned}
\end{equation}

\begin{theorem}\label{prop:mp-recovery}
 Suppose that~\eqref{eq:ose-standard} holds on $\mathscr U=\range(H)$. Under the rounding model above, assume $\epsilon>0$ and $d_s\le\epsilon/2$. Let
\begin{equation}\label{eq:mp-theta}
 \theta=4(1+\eta)\left(\frac{3d_x}{\sqrt{1-\eta}}+d_x^2+d_g\right)<1.
\end{equation}
If the following scale separation holds:
\begin{equation}\label{eq:mp-scale}
 \frac{(\epsilon+d_s)^2}{1-\eta}
 \le\frac{(\sigma_r(\widetilde Y)-d_s)^2}{2(1+\eta)},
\end{equation}
then the   projector $\widetilde P_r=\widetilde V_r\widetilde V_r\T$ satisfies
\begin{equation}\label{eq:mp-refined-interface}
 \norm{H-H\widetilde P_r}_F^2
 \le\frac{1+\theta}{1-\theta}\tau_{r+1}^2(H)
 +2r\left[
 \eta\frac{\epsilon+d_s}{\sqrt{1-\eta}}
 +3\sqrt{1+\eta}d_s
 \right]^2.
\end{equation}
{\color{black}
The bounds $d_s,d_x,d_g$ depend explicitly on $u_\ell$ through~\eqref{eq:mp-rounding-bounds}. If furthermore $S$ has bounded entries and $\epsilon=\Theta\!\left(\sqrt{u_\ell}\norm{H}_F\right)$, \eqref{eq:mp-refined-interface} then simplifies to
\begin{equation}\label{eq:mp-asymptotic-accuracy}
 \norm{H-H\widetilde P_r}_F^2
 \le\bigl(1+O(\sqrt{u_\ell})\bigr)\tau_{r+1}^2(H)+O(u_\ell).
\end{equation}
}
\end{theorem}

  $\theta$ in~\eqref{eq:mp-theta} bounds the perturbation from preconditioning multiplication and Gram formation relative to the exact Gram matrix. Requiring $\theta<1$ keeps the computed Gram matrix positive definite for Cholesky recovery. This is usually easy to satisfy for moderate $\eta$ and small accumulated rounding errors. Since~\eqref{eq:mp-rounding-bounds} uses worst-case estimates, \eqref{eq:mp-theta} is a conservative sufficient condition.

The condition $d_s\le\epsilon/2$ places the threshold above
the sketch rounding-error scale.
Condition~\eqref{eq:mp-scale}, equivalently
$\sigma_r(\widetilde Y)\ge
d_s+\sqrt{2(1+\eta)/(1-\eta)}(\epsilon+d_s)$,
requires the target singular value to be sufficiently above
the threshold. It can be   satisfied when $(1+\eta)/(1-\eta)$ is moderate and the target singular value is well above the threshold. If the target reaches the numerical truncation scale, one may reduce $r$.

In exact arithmetic, $d_s=d_x=d_g=0$ and $\theta=0$, so~\eqref{eq:mp-refined-interface} reduces to the error bound
in Theorem~\ref{thm:rp-refined}. With low-precision recovery, \eqref{eq:mp-asymptotic-accuracy} shows that the best rank-$r$ tail energy is multiplied by only $1+O(\sqrt{u_\ell})$, with an additional $O(u_\ell)$   error from truncation and sketch perturbations. This bound characterizes the effects of low-precision computations on the final approximation error.

\begin{proof}[Proof of Theorem~\ref{prop:mp-recovery}]
Write $\widetilde P_\epsilon=\widetilde V_\epsilon\widetilde V_\epsilon\T$. As in the proof of Proposition \ref{thm:rp}, $\widetilde P_\epsilon\widetilde P_r=\widetilde P_r$ and $H-H\widetilde P_r=H(I-\widetilde P_\epsilon)+H\widetilde P_\epsilon(I-\widetilde P_r)$, and
\begin{equation}\label{eq:mp-residual-split}
 \norm{H-H\widetilde P_r}_F^2
 =\norm{H(I-\widetilde P_\epsilon)}_F^2
 +\norm{H\widetilde P_\epsilon(I-\widetilde P_r)}_F^2.
\end{equation}
The first term is the energy discarded; the second is the recovery residual. Our proof idea is to first use  a   Gram comparison to bound the second term by the best rank-$r$ residual, plus a low-precision   term. We then use Lemma~\ref{lem:retained-subspace-error} and a noisy-sketch coupling estimate to bound the exact truncation error.

We first bound the singular values of $X_\circ$. Let $\widetilde U_\epsilon\in\R^{s\times t}$ contain the left
singular vectors of $\widetilde Y$ corresponding to singular
values greater than $\epsilon$, so that  $\widetilde Y\widetilde V_\epsilon=\widetilde U_\epsilon\widetilde\Sigma_\epsilon$. Using $X_\circ=H\widetilde T_\epsilon$, $SH=\widetilde Y-E_s$, and $\widetilde T_\epsilon=\widetilde V_\epsilon\widetilde\Sigma_\epsilon^{-1}$ gives
\[
 \begin{aligned}
 SX_\circ
 &=(\widetilde Y-E_s)\widetilde T_\epsilon
  =\widetilde Y\widetilde V_\epsilon\widetilde\Sigma_\epsilon^{-1}
   -E_s\widetilde T_\epsilon
 =\widetilde U_\epsilon-E_s\widetilde T_\epsilon.
 \end{aligned}
\]
Truncation in Algorithm~\ref{alg:mixed-rp} gives $\sigma_t(\widetilde Y)>\epsilon$, which together with  $\norm{\widetilde T_\epsilon}_2=1/\sigma_t(\widetilde Y)$, $d_s\le\epsilon/2$ and~\eqref{eq:mp-rounding-bounds} gives $\norm{E_s\widetilde T_\epsilon}_2\le1/2$. Hence, for $z\in\R^t$, $\norm{E_s\widetilde T_\epsilon z}_2\le\norm{z}_2/2$, while $\norm{\widetilde U_\epsilon z}_2=\norm{z}_2$. Applying the reverse triangle and triangle inequalities to $SX_\circ z=\widetilde U_\epsilon z-E_s\widetilde T_\epsilon z$ gives
\[
 \frac12\norm{z}_2\le\norm{SX_\circ z}_2\le\frac32\norm{z}_2.
\]
 Since $X_\circ z\in\range(H)$, \eqref{eq:ose-standard} gives $\sqrt{1-\eta}\norm{X_\circ z}_2\le\norm{SX_\circ z}_2\le\sqrt{1+\eta}\norm{X_\circ z}_2$. Combining the two bounds yields
\begin{equation}\label{eq:mp-conditioning}
 \frac{1}{2\sqrt{1+\eta}}\le\sigma_{\min}(X_\circ)
 \le\sigma_{\max}(X_\circ)\le\frac{3}{2\sqrt{1-\eta}}.
\end{equation}
 Next we compare the computed Gram matrix with $X_\circ\T X_\circ$.  \eqref{eq:mp-errors} gives
\[
 \widetilde G_\epsilon-X_\circ\T X_\circ
 =X_\circ\T E_x+E_x\T X_\circ+E_x\T E_x+E_g.
\]
 Use $\norm{X_\circ}_2\le3/(2\sqrt{1-\eta})$ from~\eqref{eq:mp-conditioning}, the   bounds~\eqref{eq:mp-rounding-bounds}, and $\theta$ in~\eqref{eq:mp-theta} to obtain
\[
 \begin{aligned}
 \norm{\widetilde G_\epsilon-X_\circ\T X_\circ}_2
 &\le2\norm{X_\circ}_2\norm{E_x}_2+\norm{E_x}_2^2+\norm{E_g}_2\\
 &\le\frac{3d_x}{\sqrt{1-\eta}}+d_x^2+d_g=\frac{\theta}{4(1+\eta)}.
 \end{aligned}
\]

 On the other hand, the lower bound in~\eqref{eq:mp-conditioning} gives $X_\circ\T X_\circ\succeq I/(4(1+\eta))$.
 Combining these relations gives
\[
 (1-\theta)X_\circ\T X_\circ
 \preceq\widetilde G_\epsilon
 \preceq(1+\theta)X_\circ\T X_\circ.
\]
In particular, $\widetilde G_\epsilon$ is positive definite. Multiplying on the left by $\widetilde V_\epsilon\widetilde\Sigma_\epsilon$ and on the right by $\widetilde\Sigma_\epsilon\widetilde V_\epsilon\T$ preserves these inequalities. Since $X_\circ=H\widetilde V_\epsilon\widetilde\Sigma_\epsilon^{-1}$, we have $X_\circ\widetilde\Sigma_\epsilon\widetilde V_\epsilon\T=H\widetilde V_\epsilon\widetilde V_\epsilon\T=H\widetilde P_\epsilon$. Also, $\widetilde B_\epsilon=\widetilde R_\epsilon\widetilde\Sigma_\epsilon\widetilde V_\epsilon\T$ in Algorithm~\ref{alg:mixed-rp}, with $\widetilde G_\epsilon=\widetilde R_\epsilon\T\widetilde R_\epsilon$. These  give  $\widetilde V_\epsilon\widetilde\Sigma_\epsilon\widetilde G_\epsilon\widetilde\Sigma_\epsilon\widetilde V_\epsilon\T=(\widetilde R_\epsilon\widetilde\Sigma_\epsilon\widetilde V_\epsilon\T)\T(\widetilde R_\epsilon\widetilde\Sigma_\epsilon\widetilde V_\epsilon\T)=\widetilde B_\epsilon\T\widetilde B_\epsilon$. Thus
\begin{align}\label{eq:mp-transformed-inequality}
 (1-\theta)(H\widetilde P_\epsilon)\T(H\widetilde P_\epsilon)
 \preceq\widetilde B_\epsilon\T\widetilde B_\epsilon
 \preceq(1+\theta)(H\widetilde P_\epsilon)\T(H\widetilde P_\epsilon).
\end{align}
Multiply the left inequality by $(I-\widetilde P_r)\T$ on the left and $I-\widetilde P_r$ on the right, then take traces to obtain
\[
\begin{aligned}
 (1-\theta)\norm{H\widetilde P_\epsilon(I-\widetilde P_r)}_F^2
 &=(1-\theta)\operatorname{tr}\!\left((I-\widetilde P_r)\T
 (H\widetilde P_\epsilon)\T(H\widetilde P_\epsilon)(I-\widetilde P_r)\right)\\
 &\le\operatorname{tr}\!\left((I-\widetilde P_r)\T
 \widetilde B_\epsilon\T\widetilde B_\epsilon(I-\widetilde P_r)\right)
 =\norm{\widetilde B_\epsilon(I-\widetilde P_r)}_F^2.
\end{aligned}
\]
To proceed,
let $P_*$ be a best rank-$r$ right orthogonal projector for $H\widetilde P_\epsilon$, i.e., $H\widetilde P_\epsilon P_*=[H\widetilde P_\epsilon]_r$.  Applying the right inequality of~\eqref{eq:mp-transformed-inequality} thus gives
\begin{small}
\begin{align*}
 \norm{\widetilde B_\epsilon(I-P_*)}_F^2
 &=\operatorname{tr}\!\left((I-P_*)\T\widetilde B_\epsilon\T
 \widetilde B_\epsilon(I-P_*)\right)\\
 &\le(1+\theta)\operatorname{tr}\!\left((I-P_*)\T
 (H\widetilde P_\epsilon)\T(H\widetilde P_\epsilon)(I-P_*)\right)\\
 &=(1+\theta)\norm{H\widetilde P_\epsilon(I-P_*)}_F^2
 =(1+\theta)\norm{H\widetilde P_\epsilon-[H\widetilde P_\epsilon]_r}_F^2,
\end{align*}
\end{small}
The last equality uses the definition of  $P_*$.
On the other hand, optimality of $\widetilde P_r$ for $\widetilde B_\epsilon$ gives $\norm{\widetilde B_\epsilon(I-\widetilde P_r)}_F^2\le\norm{\widetilde B_\epsilon(I-P_*)}_F^2$. This   in connection with the above two relations gives
\[
\begin{aligned}
 \norm{H\widetilde P_\epsilon(I-\widetilde P_r)}_F^2
 & \le\frac{1+\theta}{1-\theta}
 \norm{H\widetilde P_\epsilon-[H\widetilde P_\epsilon]_r}_F^2\\
 &=\norm{H\widetilde P_\epsilon-[H\widetilde P_\epsilon]_r}_F^2
 +\frac{2\theta}{1-\theta}
 \norm{H\widetilde P_\epsilon-[H\widetilde P_\epsilon]_r}_F^2,
\end{aligned}
\]
This bounds the second term of~\eqref{eq:mp-residual-split}. We then consider the first term of~\eqref{eq:mp-residual-split}.
Combine it with the first term above. Then orthogonality yields
\[
 \norm{H(I-\widetilde P_\epsilon)}_F^2
 +\norm{H\widetilde P_\epsilon-[H\widetilde P_\epsilon]_r}_F^2
 =\norm{H-[H\widetilde V_\epsilon]_r\widetilde V_\epsilon\T}_F^2,
\]
where we have used $H\widetilde P_\epsilon=H\widetilde V_\epsilon\widetilde V_\epsilon\T$ and $[H\widetilde P_\epsilon]_r=[H\widetilde V_\epsilon]_r\widetilde V_\epsilon\T$; here the singular vectors are chosen consistently. Substitution into~\eqref{eq:mp-residual-split} therefore gives
\begin{align}\label{eq:mp-error-reduction}
 \norm{H-H\widetilde P_r}_F^2
 &\le\norm{H-[H\widetilde V_\epsilon]_r\widetilde V_\epsilon\T}_F^2
 +\frac{2\theta}{1-\theta}
 \norm{H\widetilde P_\epsilon-[H\widetilde P_\epsilon]_r}_F^2\nonumber\\
 &\leq \norm{H-[H\widetilde V_\epsilon]_r\widetilde V_\epsilon\T}_F^2
 +\frac{2\theta}{1-\theta}
 \tau_{r+1}^2(H),
\end{align}
where   the second inequality  uses $\norm{H\widetilde P_\epsilon-[H\widetilde P_\epsilon]_r}_F=\tau_{r+1}(H\widetilde P_{\epsilon})$ and $\tau_{r+1}(HP) \leq \tau_{r+1}(H)$ for any orthogonal projector $P$.

We will next bound $\|{H-[H\widetilde V_\epsilon]_r\widetilde V_\epsilon\T}\|_F^2$ using Lemma~\ref{lem:retained-subspace-error}. Complete
$\widetilde V_\epsilon$ to an orthogonal matrix
$[\widetilde V_\epsilon,\widetilde V_\perp]$
using the remaining right singular vectors of $\widetilde Y$.
We first verify the assumptions of the lemma below.

  Note that $\sigma_{\min}(X_\circ)>0$ by~\eqref{eq:mp-conditioning}. Since $\widetilde\Sigma_\epsilon$ is invertible, $H\widetilde V_\epsilon=X_\circ\widetilde\Sigma_\epsilon$ has full column rank, as required with $V=\widetilde V_\epsilon$ in Lemma~\ref{lem:retained-subspace-error}. For the scale condition, truncation gives $\norm{\widetilde Y\widetilde V_\perp}_2\le\epsilon$, so
 \begin{align}\label{eq:mp-discarded-subspace}
 \sqrt{1-\eta}\norm{H\widetilde V_\perp}_2
 &\le\norm{SH\widetilde V_\perp}_2
  \le\norm{\widetilde Y\widetilde V_\perp}_2
      +\norm{E_s\widetilde V_\perp}_2
 \le\epsilon+d_s,
 \end{align}
where the first inequality follows from the lower bound in~\eqref{eq:ose-standard} as $\range(H\widetilde V_\perp)\subseteq\mathscr U$; the second uses $SH=\widetilde Y-E_s$, and the last uses $\norm{E_s}_2\le d_s$.
On the other hand, the singular values of $\widetilde Y\widetilde V_\epsilon$ are the leading $t$ singular values of $\widetilde Y$. From $SH\widetilde V_\epsilon=\widetilde Y\widetilde V_\epsilon-E_s\widetilde V_\epsilon$ and the singular-value perturbation bound, for $1\le j\le t$,
\begin{equation}\label{eq:mp-retained-sketch-sv}
 \sigma_j(SH\widetilde V_\epsilon)
 \ge\sigma_j(\widetilde Y\widetilde V_\epsilon)
      -\norm{E_s\widetilde V_\epsilon}_2
 \ge\sigma_j(\widetilde Y)-d_s.
\end{equation}
The upper bound in~\eqref{eq:ose-standard} on $\range(H\widetilde V_\epsilon)$ gives $(SH\widetilde V_\epsilon)\T(SH\widetilde V_\epsilon)\preceq(1+\eta)(H\widetilde V_\epsilon)\T(H\widetilde V_\epsilon)$, which together with~\eqref{eq:mp-retained-sketch-sv} yields
\begin{equation}\label{eq:mp-retained-original-sv}
 \sigma_j(H\widetilde V_\epsilon)
 \ge \frac{1}{\sqrt{1+\eta}}{\sigma_j(SH\widetilde V_\epsilon)}
 \ge  \frac{1}{\sqrt{1+\eta}}({\sigma_j(\widetilde Y)-d_s}) ,\qquad 1\le j\le t.
\end{equation}
 Taking $j=r$ and using~\eqref{eq:mp-scale} and~\eqref{eq:mp-discarded-subspace}, we get
\[
 \norm{H\widetilde V_\perp}_2^2
 \le\frac{(\epsilon+d_s)^2}{1-\eta}
 \le\frac{(\sigma_r(\widetilde Y)-d_s)^2}{2(1+\eta)}
 \le\frac12\sigma_r^2(H\widetilde V_\epsilon).
\]
Therefore,
the scale condition of Lemma~\ref{lem:retained-subspace-error} now holds.

Now, take a thin QR factorization $H\widetilde V_\epsilon=\widetilde Q C$. Lemma  \ref{lem:retained-subspace-error} applies with $V=\widetilde V_\epsilon$ and the orthonormal basis $Q=\widetilde Q$, giving

\begin{equation}\label{eq:mp-ideal-truncation}
 \norm{H-[H\widetilde V_\epsilon]_r\widetilde V_\epsilon\T}_F^2
 \le\tau_{r+1}^2(H)+2r\norm{\widetilde Q\T H\widetilde V_\perp}_2^2.
\end{equation}
To bound   $\norm{\widetilde Q\T H\widetilde V_\perp}_2$, take $j=t$ in~\eqref{eq:mp-retained-original-sv}. With $\sigma_t(\widetilde Y)>\epsilon$ and $d_s\le\epsilon/2$, this gives
\begin{align}\label{eq:mp-inverse-bound}
 \norm{C^{-1}}_2
 =\frac{1}{\sigma_t(H\widetilde V_\epsilon)}
 \le\frac{\sqrt{1+\eta}}{\sigma_t(\widetilde Y)-d_s}
 \le\frac{\sqrt{1+\eta}}{\epsilon-d_s}.
\end{align}
From $\widetilde Q=H\widetilde V_\epsilon C^{-1}$ and $\widetilde Y=SH+E_s$ in~\eqref{eq:mp-errors},
\[
 S\widetilde Q+E_s\widetilde V_\epsilon C^{-1}
 =(SH+E_s)\widetilde V_\epsilon C^{-1}
 =\widetilde Y\widetilde V_\epsilon C^{-1}.
\]
As $\widetilde V_\epsilon$ and $\widetilde V_\perp$ come from
the same SVD of $\widetilde Y$, their columns are orthogonal
eigenvectors of $\widetilde Y\T\widetilde Y$, so
$(\widetilde Y\widetilde V_\epsilon)\T
\widetilde Y\widetilde V_\perp=0$.
Combining this with the above identity leads to
\[
 (S\widetilde Q+E_s\widetilde V_\epsilon C^{-1})\T
 \widetilde Y\widetilde V_\perp
 =C^{-\mathsf T}(\widetilde Y\widetilde V_\epsilon)\T
 \widetilde Y\widetilde V_\perp=0.
\]
Rearranging this   and substituting $SH\widetilde V_\perp=\widetilde Y\widetilde V_\perp-E_s\widetilde V_\perp$ yields
\[
 \begin{aligned}
 (S\widetilde Q)\T SH\widetilde V_\perp
 &=(S\widetilde Q)\T\widetilde Y\widetilde V_\perp
   -(S\widetilde Q)\T E_s\widetilde V_\perp
 \\
  & =-(E_s\widetilde V_\epsilon C^{-1})\T\widetilde Y\widetilde V_\perp
   -(S\widetilde Q)\T E_s\widetilde V_\perp.
 \end{aligned}
\]
Using $\norm{E_s}_2\le d_s$ and $\norm{\widetilde Y\widetilde V_\perp}_2\le\epsilon$ gives
\[
 \begin{aligned}
 \norm{(S\widetilde Q)\T SH\widetilde V_\perp}_2
 &\le\norm{E_s\widetilde V_\epsilon C^{-1}}_2
       \norm{\widetilde Y\widetilde V_\perp}_2
       +\norm{S\widetilde Q}_2\norm{E_s\widetilde V_\perp}_2\\
 &\le d_s\norm{C^{-1}}_2\epsilon+\norm{S\widetilde Q}_2d_s.
 \end{aligned}
\]
Since $\widetilde Q$ is orthonormal with columns in $\range(H)$, the upper bound in~\eqref{eq:ose-standard} gives $\norm{S\widetilde Q}_2\le\sqrt{1+\eta}$. Combining the bound of $\|C^{-1}\|_2$ in \eqref{eq:mp-inverse-bound} with $\epsilon/(\epsilon-d_s)\le2$ from $d_s\le\epsilon/2$ yields
\begin{align}\label{eq:mp-sketch-coupling}
 \norm{(S\widetilde Q)\T SH\widetilde V_\perp}_2
 \le\sqrt{1+\eta}\,d_s\left(1+\frac{\epsilon}{\epsilon-d_s}\right)
 \le3\sqrt{1+\eta}\,d_s.
\end{align}

To transfer this estimate to the original space, take an orthonormal basis $Q_H$ for $\range(H)$. The embedding condition~\eqref{eq:ose-standard} is equivalent to
\[
 -\eta I
 \preceq I-Q_H\T S\T SQ_H
 \preceq\eta I,
\]
namely, $\|I-Q_H\T S\T SQ_H\|_2\leq \eta$. On the other hand,    since the columns of $\widetilde Q$ and $H\widetilde V_\perp$ lie in $\range(H)$, we have  $\widetilde Q^T=\widetilde Q^TQ_HQ_H^T,
H\widetilde V_\perp=Q_HQ_H^TH\widetilde V_\perp$.  Thus
\begin{small}
 \begin{align}\label{eq:mp-projection-coupling}
  \norm{\widetilde Q\T\left(I-S\T S\right)H\widetilde V_\perp}_2
   =\norm{\widetilde Q\T Q_H
 \left(I-Q_H\T S\T SQ_H\right)
 Q_H\T H\widetilde V_\perp}_2  \le\eta\norm{H\widetilde V_\perp}_2.
 \end{align}
\end{small}
Write $\widetilde Q\T H\widetilde V_\perp=\widetilde Q\T(I-S\T S)H\widetilde V_\perp+(S\widetilde Q)\T SH\widetilde V_\perp$. Apply~\eqref{eq:mp-projection-coupling} and~\eqref{eq:mp-sketch-coupling} to the two terms and use~\eqref{eq:mp-discarded-subspace} to obtain
\begin{equation}\label{eq:mp-noisy-coupling}
 \begin{aligned}
 \norm{\widetilde Q\T H\widetilde V_\perp}_2
  \le\eta\norm{H\widetilde V_\perp}_2
 +\norm{(S\widetilde Q)\T SH\widetilde V_\perp}_2 \le\eta\frac{\epsilon+d_s}{\sqrt{1-\eta}}
 +3\sqrt{1+\eta}d_s.
 \end{aligned}
\end{equation}
Combining~\eqref{eq:mp-noisy-coupling}, \eqref{eq:mp-ideal-truncation}, and~\eqref{eq:mp-error-reduction} proves~\eqref{eq:mp-refined-interface}.

 Finally, we prove~\eqref{eq:mp-asymptotic-accuracy}. Suppose $|S_{ij}|\le K$.   \eqref{eq:mp-rounding-bounds} then gives
\[
 d_s=\gamma_{M+2}^{(\ell)}\norm{|S||H|}_2
 \le\gamma_{M+2}^{(\ell)}K\sqrt{sM}\norm{H}_F
 =O(u_\ell)\norm{H}_F.
\]
It follows from the   truncation rule and~\eqref{eq:mp-conditioning} that
\[
 \begin{aligned}
 d_x&\le\gamma_{n+2}^{(\ell)}\norm{H}_F\norm{\widetilde T_\epsilon}_F
 \le\gamma_{n+2}^{(\ell)}\frac{\sqrt t\norm{H}_F}{\epsilon},\\
 d_g&\le\gamma_{M+2}^{(\ell)}\norm{\widetilde X_\epsilon}_F^2
 \le t\gamma_{M+2}^{(\ell)}
       \left(\frac{3}{2\sqrt{1-\eta}}+d_x\right)^2.
 \end{aligned}
\]
The above uses $\sigma_t(\widetilde Y)>\epsilon$, $\widetilde X_\epsilon=X_\circ+E_x$, and $t\le n$. Since $\gamma_j^{(\ell)}=O(u_\ell)$, substituting $\epsilon=\Theta\!\left(\sqrt{u_\ell}\norm{H}_F\right)$ into~\eqref{eq:mp-rounding-bounds} and~\eqref{eq:mp-theta} yields successively
\[
 d_s=O(u_\ell),\qquad d_x=O(\sqrt{u_\ell}),\qquad
 d_g=O(u_\ell),\qquad\theta=O(\sqrt{u_\ell}).
\]
 Thus $(1+\theta)/(1-\theta)=1+O(\sqrt{u_\ell})$, and the second term of~\eqref{eq:mp-refined-interface} is $O((\epsilon+d_s)^2)=O(u_\ell)$. Substituting these  into~\eqref{eq:mp-refined-interface} proves~\eqref{eq:mp-asymptotic-accuracy}.
\end{proof}

\section{Numerical Experiments}\label{sec:numerical-experiments}

Results are plotted as medians in this section, while complete numerical results are tabulated as means with sample standard deviations in the supplementary material.

We test the proposed algorithms in this section; the mixed-precision variants use the mixed single and double precision.   All   experiments employ MATLAB R2025a on an Apple M3 Ultra with 512\,GB RAM.

\subsection{Experimental setup}\label{sec:experimental-setup}

\subsubsection{Compared algorithms} \label{sec:experiment-comparison}
Classical HOSVD and STHOSVD serve as the primary references,
since our aim is to accelerate these decompositions while
preserving their accuracy. We also compare with randomized
Tucker approximation methods to assess whether the proposed
methods can retain classical-level accuracy at a comparable
computational cost (recall from
Sect.~\ref{sec:rangefinding-recovery} that randomization
serves a different role in our methods). The comparison methods are listed below:

  MSK~\cite{MSK2020} (randomized SVD), MLB~\cite{MLB2024} (two-stage Tucker), KRP~\cite{SVB2026} (Khatri--Rao projections), Shift~\cite{CWWY2025} (adaptive shifts), Stable-MLN~\cite{BR2024} (stabilized multilinear Nystr\"om), and RTSMS~\cite{HN2023} (single-mode sketching) \footnote{We adapt the authors' codes for \href{https://github.com/rlminste/rand_tensor_decompositions}{MSK}, \href{https://github.com/bhisham123/KRP-Random-Projections}{KRP}, \href{https://github.com/bhashemi/rtsms}{RTSMS}, and \href{https://github.com/alb95/MLN}{Stable-MLN}. MLB and Shift are implemented from their published algorithms. We also implement the nonsequential RTSMS variant and sequential Stable-MLN~\cite{BH2024}. }.

\paragraph{Parameters.} We set over sampling $p=10$ for MSK, MLB and Shift, matching the setting recommended for Shift~\cite{CWWY2025} and exceeding those used in the original MSK and MLB experiments~\cite{MSK2020,MLB2024} to improve approximation accuracy. Shift uses one shifted power iteration. KRP uses no oversampling~\cite{SVB2026}. Stable-MLN uses $p=\min\{n_j-r_j,\max(2,\lceil r_j/2\rceil)\}$~\cite{BR2024}. RTSMS uses $p=\min\{n_j,\operatorname{round}(1.5r_j)\}-r_j$ and auxiliary least-squares oversampling factor $k=4$~\cite{HN2023}. All other parameters retain the original settings.

For the proposed methods, the sketch dimension is $s_j=2n_j$, with upward rounding for KronSS. We use the numerical threshold
$\epsilon_j=2\varepsilon_{\mathrm{mach}}\|Y_j\|_2$, where $Y_j$ is the sketch of the unfolding processed at step $j$ and
  $\varepsilon_{\mathrm{mach}}$ is MATLAB's machine epsilon for double precision in RP and single precision in MP-RP.

\subsubsection{Test matrices}\label{sec:experiment-test-matrices}
Following the practical recommendation of Martinsson and Tropp~\cite[Section~9.2]{MT2020}, we use sparse sign (SS) sketches with \(\zeta=\min\{s,8\}\) randomly placed nonzeros per column, each taking the value \(\pm1/\sqrt{\zeta}\). Such sketches offer a computationally inexpensive alternative to Gaussian random projections. Our SS routine is   from \href{https://github.com/eepperly/Iterative-Sketching-Is-Stable}{Epperly's code}. KronSS uses Kronecker products of small SS matrices~\cite{MLB2024}, allowing reuse of factors and intermediate contractions to further reduce computation.  Comparison of Gaussian, SS and KronSS sketches is given in the supplementary material.

\subsubsection{Data} \label{sec:experiment-data}
\paragraph{Synthetic data.} We construct tensors with prescribed unfolding spectra as
\begin{equation}\label{eq:experiment-cp-data}
 \mathcal A=\sum_{k=1}^{K}\lambda_k\,
 u_k^{(1)}\circ\cdots\circ u_k^{(d)},
 ~
 U^{(j)}=[u_1^{(j)},\ldots,u_K^{(j)}],~
 (U^{(j)})\T U^{(j)}=I_K,
\end{equation}
where $\circ$ denotes the outer product. The factors are obtained by orthogonalizing independent Gaussian matrices. The tests cover tensors of orders 3 through 6,
with varying dimensions, spectral decay and mode shapes.
Specific settings are given in the corresponding experiments.
\begin{enumerate}
 \item \emph{Varying dimensions.} We set $n_1=\cdots=n_d=n$, with $K=n$ and target ranks $r_j=\operatorname{round}(0.2n)$.
 $\lambda_k\propto\exp[-5.98(k-1)/(K-1)]$.
 \item \emph{Varying spectral decay.} We fix the tensor size at $200\times200\times200\times200$ and $K=200$, and take $r_j=r$ with $r=10,20,\ldots,100$. The three spectra are
 \begingroup
 \setlength{\jot}{0pt}
 \begin{equation}\label{eq:experiment-spectra}
  \begin{aligned}
   \text{Slow:}\quad &\lambda_k\propto k^{-1},\\
   \text{Medium:}\quad &\lambda_k\propto e^{-5.98(k-1)/(K-1)},\\
   \text{Fast:}\quad &\lambda_k\propto e^{-17.94(k-1)/(K-1)}.
  \end{aligned}
 \end{equation}
 \endgroup
 \item \emph{Varying shapes at fixed volume.} We use ten fourth-order shapes with $\prod_{j=1}^4n_j=200^4$, from $(200,200,200,200)$ to $(100,100,400,400)$. All use $K=100$, the medium spectrum in~\eqref{eq:experiment-spectra}, and $r_j=20$.
\end{enumerate}

\paragraph{Real data.}
We use three tensors: DC Mall~\cite{DCMall}, a hyperspectral image of size
$1280\times307\times191$; Miranda~\cite{Miranda2022}, a simulation density
subtensor of size $2048\times256\times256$; and CAVE~\cite{CAVE}, comprising
32 hyperspectral scenes stacked into a
$512\times512\times31\times32$ tensor.

\subsubsection{Performance measures}
For an approximation $\widehat{\mathcal A}$ computed from the   core and factors obtained from an algorithm, let
 $
 E= {\|\mathcal A-\widehat{\mathcal A}\|_F}
 $ and $T$ the computation time; in particular,   $E_{\mathrm{classic}}$ and $T_{\mathrm{classic}}$ denote those of the classical method. We report
\begin{equation*}
 \text{Speedup}=\frac{T_{\mathrm{classic}}}{T},
 \qquad
 \text{Relative error}=\frac{E}{E_{\mathrm{classic}}}-1,
\end{equation*}
namely, the speedup and relative error over the classical HOSVD or STHOSVD. The relative error follows the   form used in~\cite[Section~7.1]{TYUC2017}.

\subsection{Numerical results for synthetic data}\label{sec:synthetic-results}

\subsubsection{Varying dimensions}
\begin{figure}[!htbp]
 \centering
 \includegraphics[width=\textwidth,height=.82\textheight,keepaspectratio]{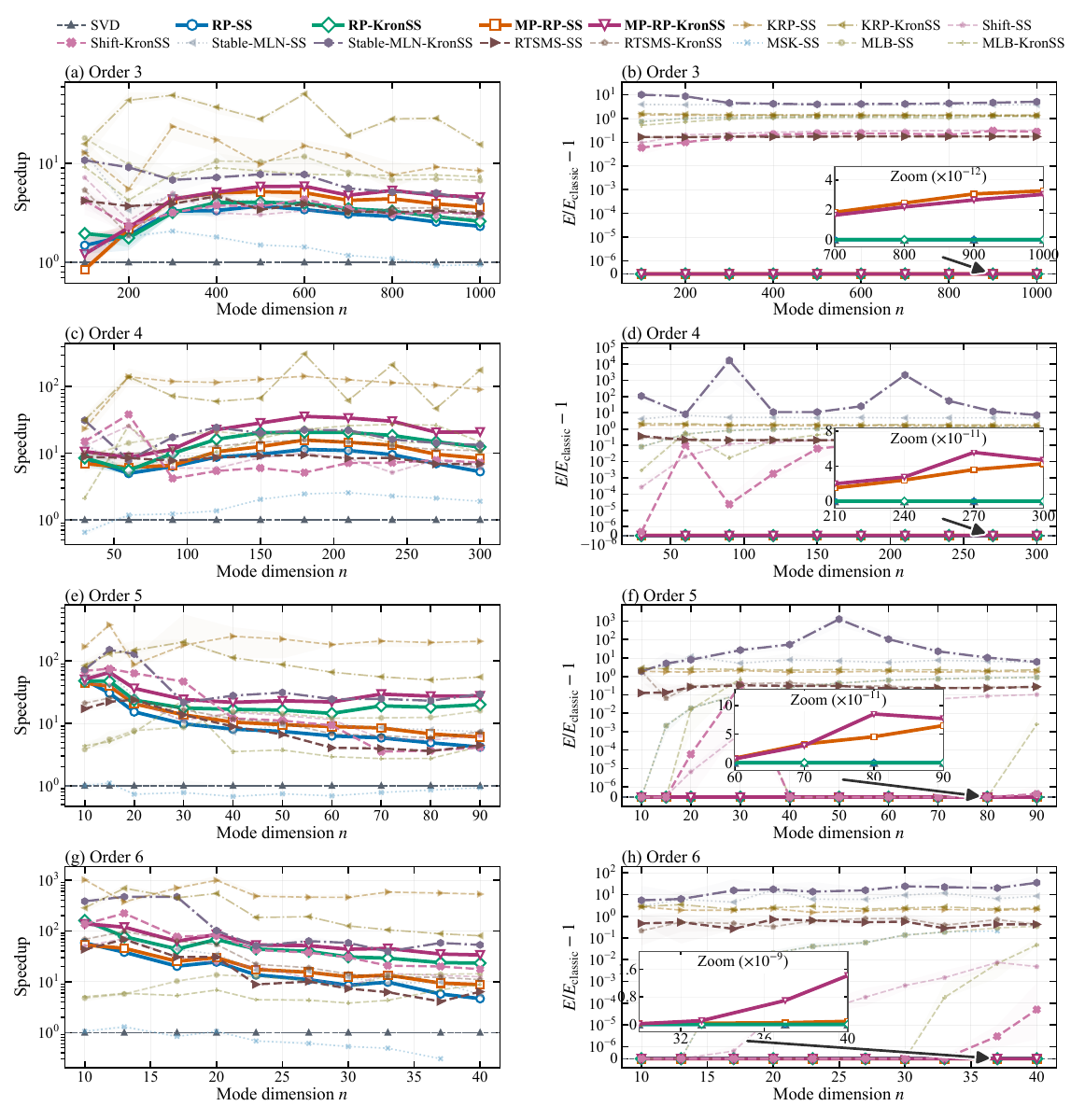}
 \caption{Varying dimensions: HOSVD comparison for orders 3 through 6. Left: speedup over classical HOSVD. Right: relative error $E/E_{\mathrm{classic}}-1$, with magnified insets. Four curves of our algorithms  are shown with thicker lines for better clarity.}
 \label{fig:synthetic-dimensions-h}
\end{figure}
\begin{figure}[!htbp]
 \centering
 \includegraphics[width=\textwidth,height=.82\textheight,keepaspectratio]{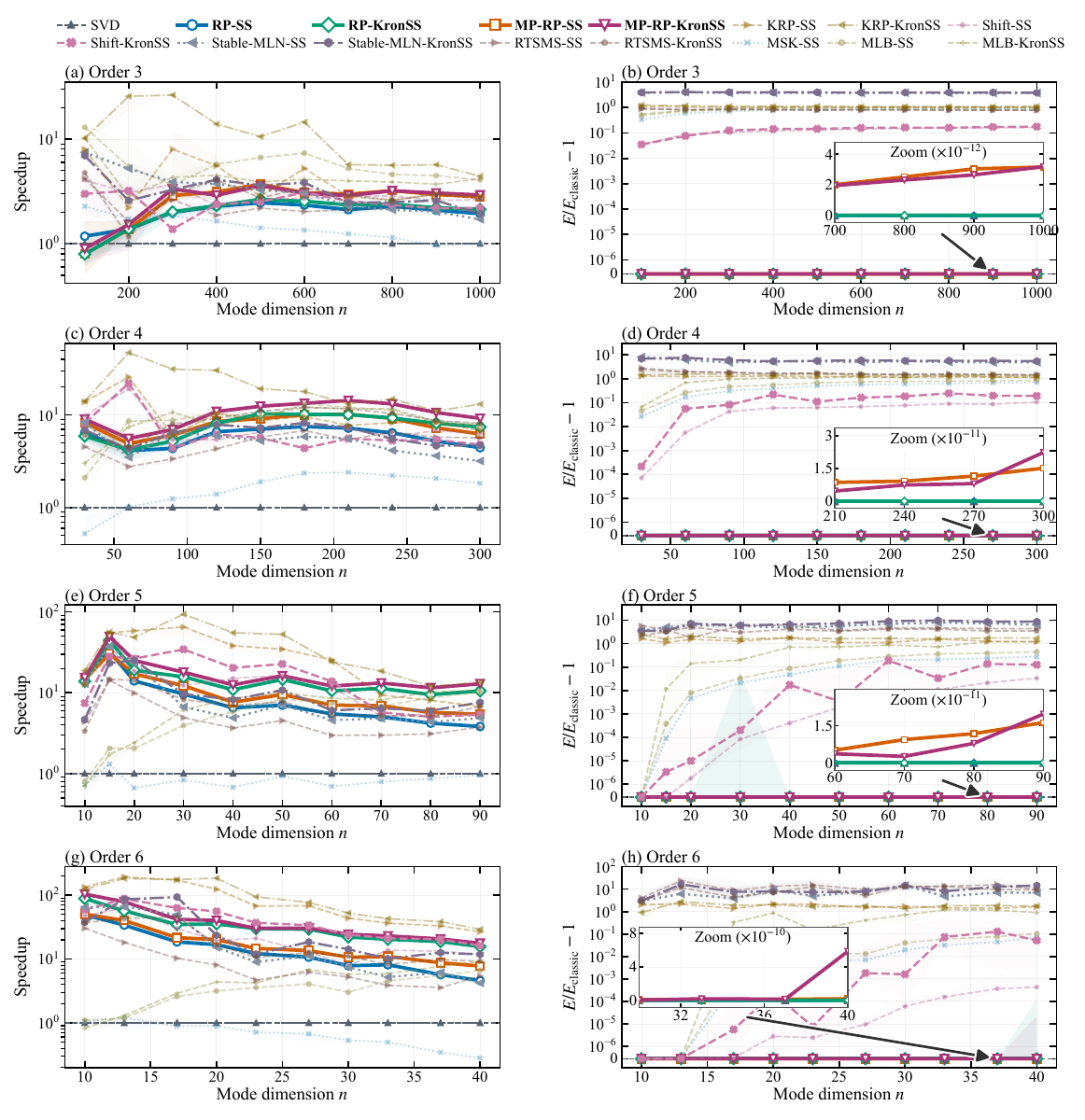}
 \caption{Varying dimensions: STHOSVD comparison.
 }
 \label{fig:synthetic-dimensions-st}
\end{figure}
Fig.~\ref{fig:synthetic-dimensions-h} (HOSVD)  and~\ref{fig:synthetic-dimensions-st} (STHOSVD) compare the methods for tensors of orders 3 through 6. Left and right panels show the speedup and relative error, respectively. The colors and markers identify the algorithms listed  on   top of Fig.~\ref{fig:synthetic-dimensions-h}.  Curves show medians, with $20$th--$80$th percentile bands where five valid runs are available. RP and MP-RP denote our ordinary and mixed-precision version; the suffixes SS and KronSS specify the test matrices.

\paragraph{Accuracy.} The right panels show that RP matches the classical recovery error to numerical precision across all tested dimensions and orders: the relative error remains indistinguishable from zero even in the magnified insets. The MP-RP curves are likewise indistinguishable from zero in the main panels. Only the magnified insets reveal minute differences, amounting to less than $0.000001\%$ in recovery error across all panels. This accuracy is consistently observed with both SS and
KronSS in the HOSVD and STHOSVD frameworks.
 Thus neither the change of test matrices nor the use of low precision introduces an appreciable loss of recovery accuracy in this experiment.

Several comparison methods, including MSK-SS, MLB-SS,
MLB-KronSS and both Shift variants, attain nearly classical
accuracy at small dimensions in both frameworks, but do not
maintain it across all tested sizes. In the HOSVD tests,
Shift-KronSS and MLB-KronSS retain this accuracy over more
dimensions. Stable-MLN-KronSS shows the largest deviations,
reflecting its single-pass recovery without revisiting the data.

These observations are consistent with the preconditioned   viewpoint discussed in Sect.~\ref{sec:rangefinding-recovery}.

\paragraph{Time.} All our algorithms    accelerate the classical decompositions substantially. MP-RP-KronSS is fastest in most settings, due to its combination of lower-cost large matrix operations  with the additional reuse enabled by KronSS; RP-KronSS is  usually the second fastest. For example, for the $13^6$ tensor, MP-RP-KronSS achieves speedups of about $118.8$ and $79.0$ over classical HOSVD and STHOSVD. Nevertheless, even RP-SS, the slowest of the four in this example, still achieves speedups of $38.2$ and $34.1$.

Compared with other algorithms, ours offer competitive computation times while retaining classical-level accuracy. KRP-SS and KRP-KronSS are faster in many settings, particularly in the HOSVD tests, but incur larger recovery errors. MSK-SS can even be slower than the classical method at larger dimensions, as it still computes an SVD of a very wide matrix despite sketching.

\subsubsection{Varying spectral decay}

\begin{figure}[!htbp]
 \centering
 \includegraphics[width=\textwidth,height=.82\textheight,keepaspectratio]{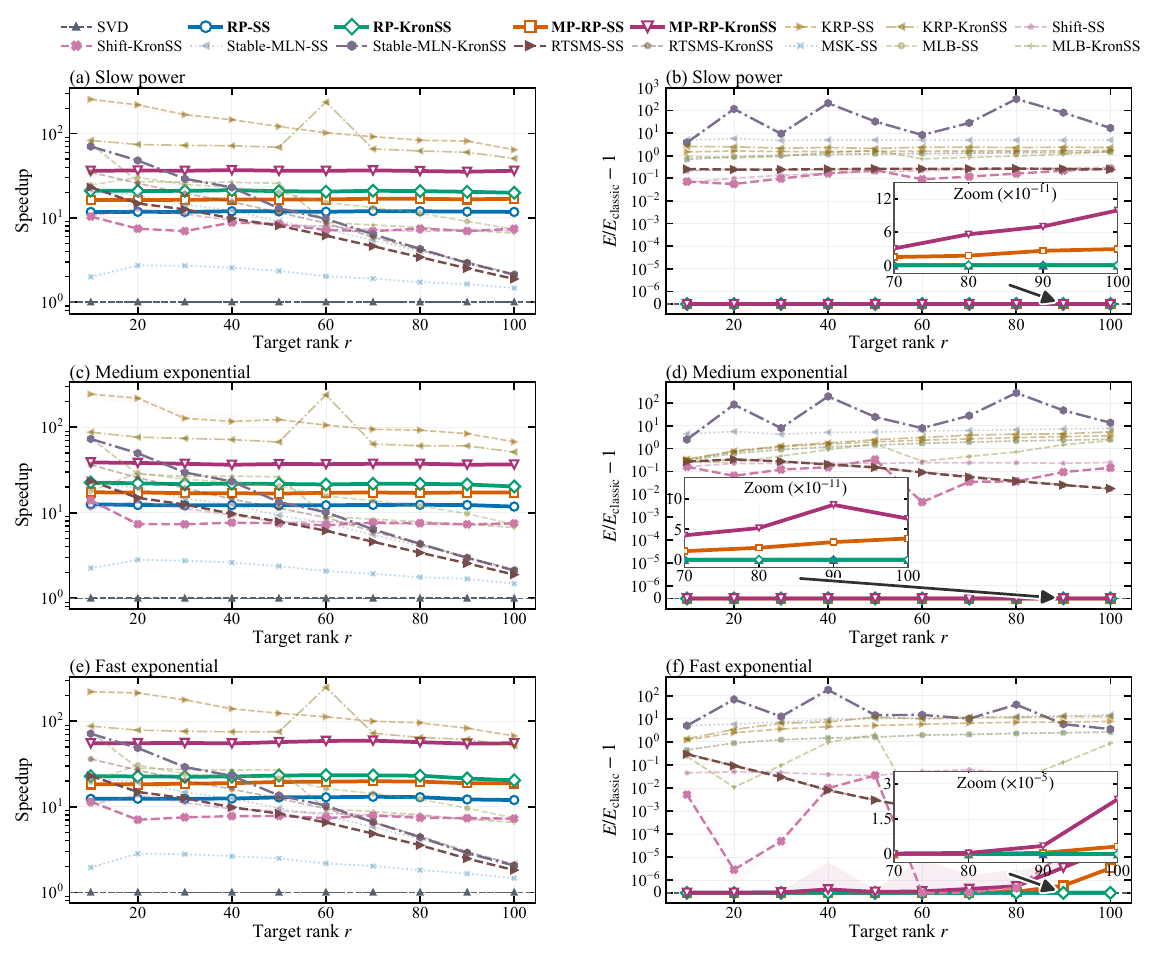}
 \caption{Varying target rank for $200^4$ tensors under slow, medium and fast spectral decay: HOSVD comparison. Left: speedup. Right: relative error and magnified insets.}
 \label{fig:synthetic-spectra-h}
\end{figure}
\begin{figure}[!htbp]
 \centering
 \includegraphics[width=\textwidth,height=.82\textheight,keepaspectratio]{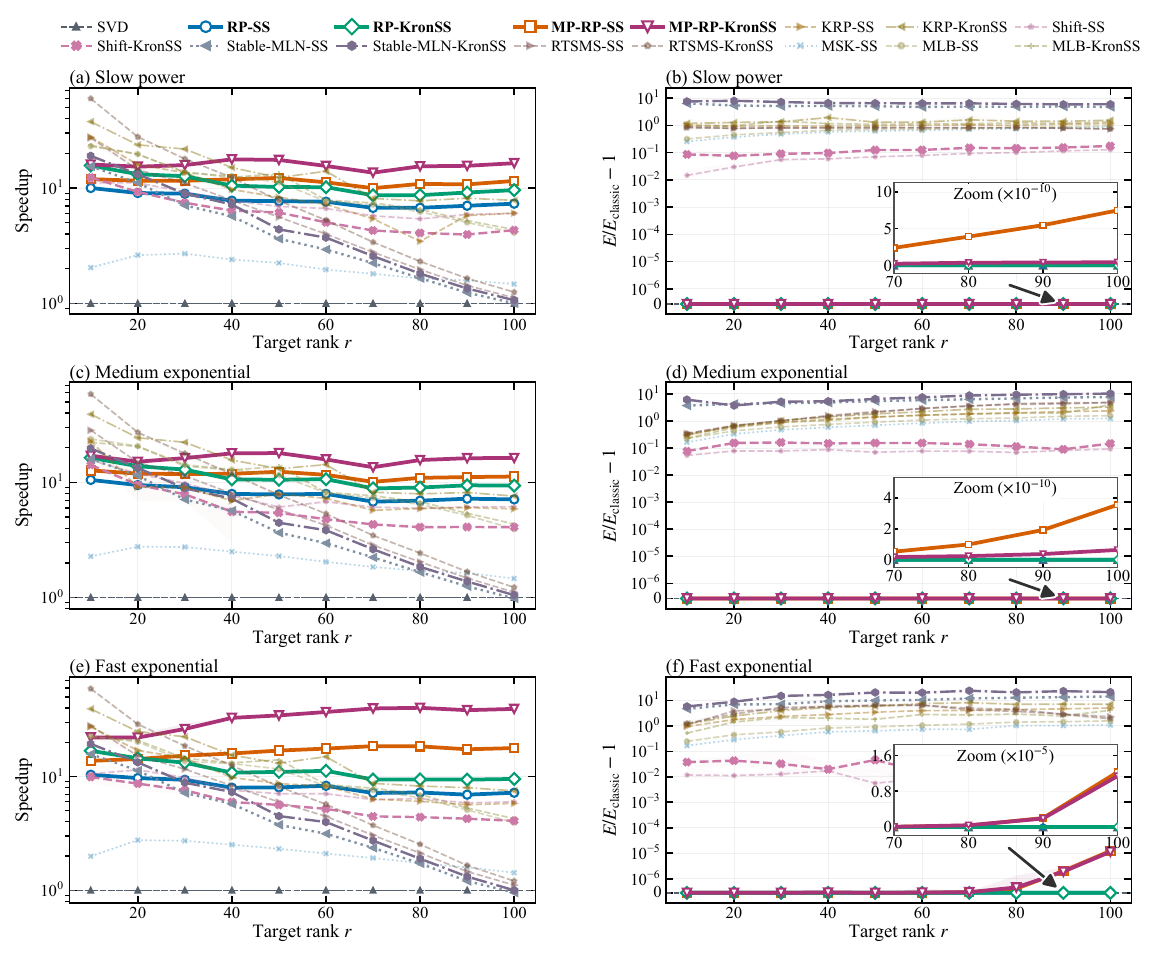}
 \caption{Varying target rank for $200^4$ tensors under slow, medium and fast spectral decay: STHOSVD comparison, with the same presentation as Fig.~\ref{fig:synthetic-spectra-h}.}
 \label{fig:synthetic-spectra-st}
\end{figure}
\paragraph{Accuracy.} Fig.~\ref{fig:synthetic-spectra-h} and~\ref{fig:synthetic-spectra-st} show that RP retains classical-level accuracy under slow, medium and fast spectral decay, with differences from zero remaining indistinguishable even in the magnified insets. MP-RP also remains close to the classical reference, with relative discrepancies below $10^{-9}$ for the medium and slow spectra. Its small deviations become visible for the fast spectrum at $r>80$,   possibly because low-precision errors account for a larger
proportion of the recovery error under fast spectral decay.  Even there, the largest relative discrepancy is about   $0.00234\%$, and the recovery errors remain below those of the   comparison methods.

\paragraph{Time.} The time panels again show substantial speedups of our algorithms over the classical decompositions across the three spectra. Particularly in the STHOSVD tests,  the speedups of several comparison methods decrease as the target rank grows, with some approaching one, while the speedups of our method remains almost unchanged. Consequently, it is interesting to see that   the proposed methods become the fastest at larger target ranks. These results again confirm substantial speedups over classical
decompositions with nearly unchanged recovery accuracy
across all three spectra.

\subsubsection{Varying shapes at fixed volume}
\begin{figure}[!htbp]
 \centering
 \includegraphics[width=\textwidth]{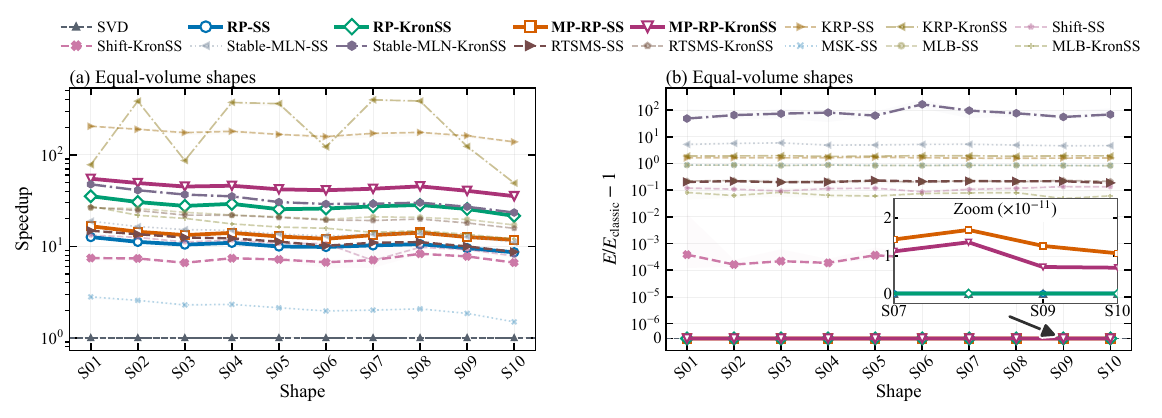}
 \par\medskip
 \includegraphics[trim=0 0 0 27bp,clip,width=\textwidth]{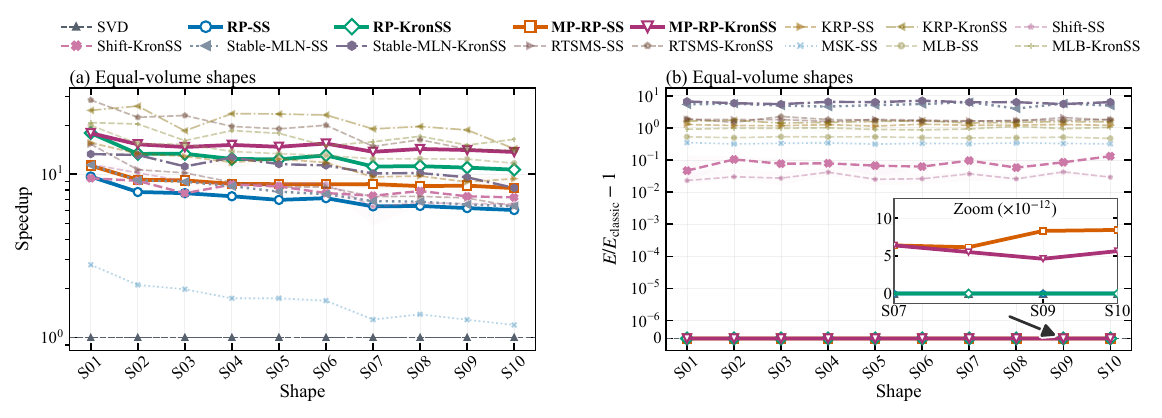}
 \caption{Fixed-volume shape comparison: HOSVD (top) and STHOSVD (bottom). Each tensor has $200^4$ entries and target ranks $(20,20,20,20)$.  }
 \label{fig:synthetic-shapes}
\end{figure}

Fig.~\ref{fig:synthetic-shapes} compares ten fixed-volume
shapes, S01--S10, with mode sizes
$(200,200,200,200)$, $(160,200,200,250)$,
$(160,160,250,250)$, $(125,200,200,320)$,
$(125,160,250,320)$, $(125,125,320,320)$,
$(100,200,250,320)$, $(100,200,200,400)$,
$(100,160,250,400)$, and $(100,100,400,400)$,
respectively. Despite the changes in tensor shape at fixed volume, the proposed methods maintain substantial speedups and classical-level recovery accuracy. The additional gains from KronSS and mixed precision also persist. For S10, KronSS reduces the   RP time relative to SS by factors of about $2.50$ and $1.77$ in the HOSVD and STHOSVD tests, respectively.

\subsection{Numerical results for real data}\label{sec:real-results}

\begin{figure}[!htbp]
 \centering
 \includegraphics[width=\textwidth,height=.82\textheight,keepaspectratio]{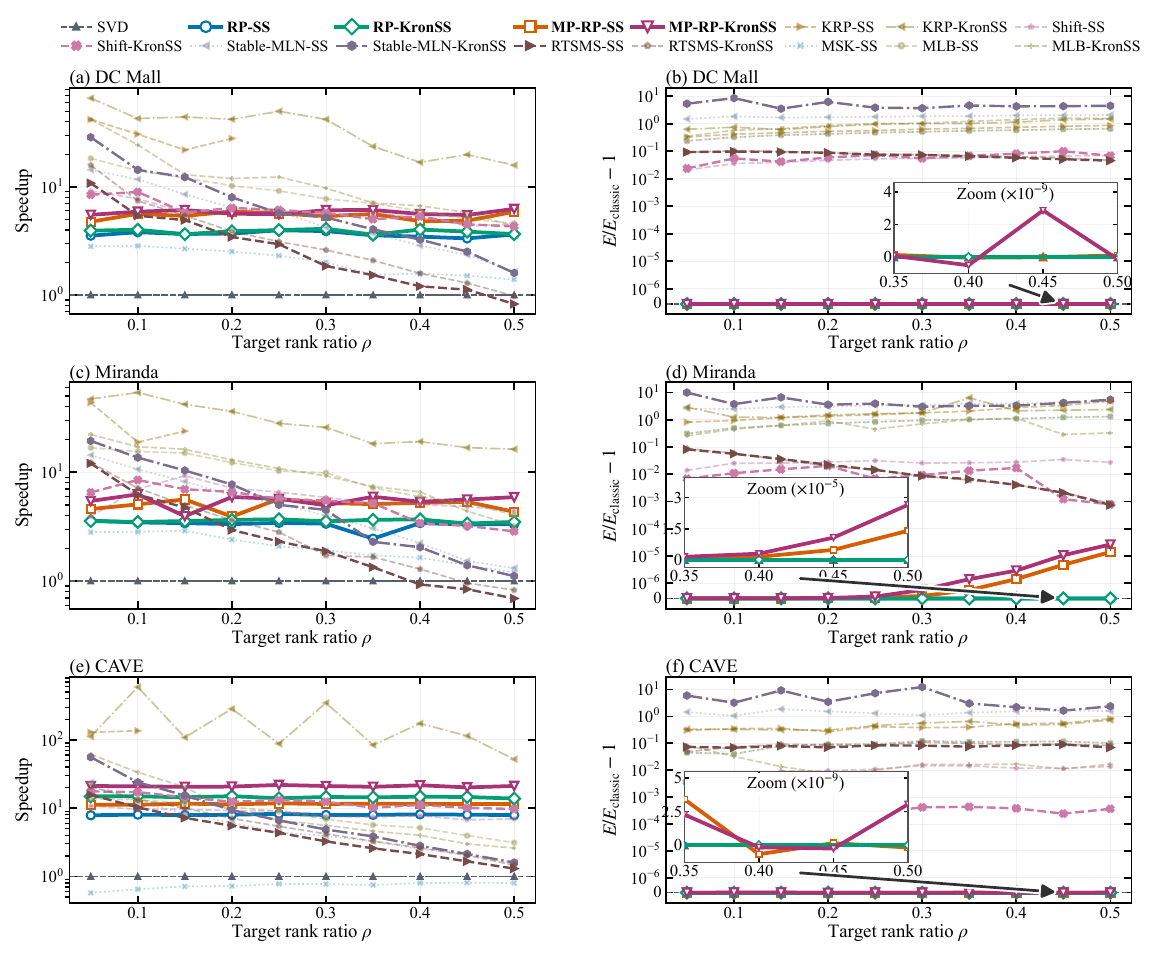}
 \caption{Real-data comparison: HOSVD tests on DC Mall, Miranda and CAVE.  }
 \label{fig:real-h}
\end{figure}

\begin{figure}[!htbp]
 \centering
 \includegraphics[width=\textwidth,height=.82\textheight,keepaspectratio]{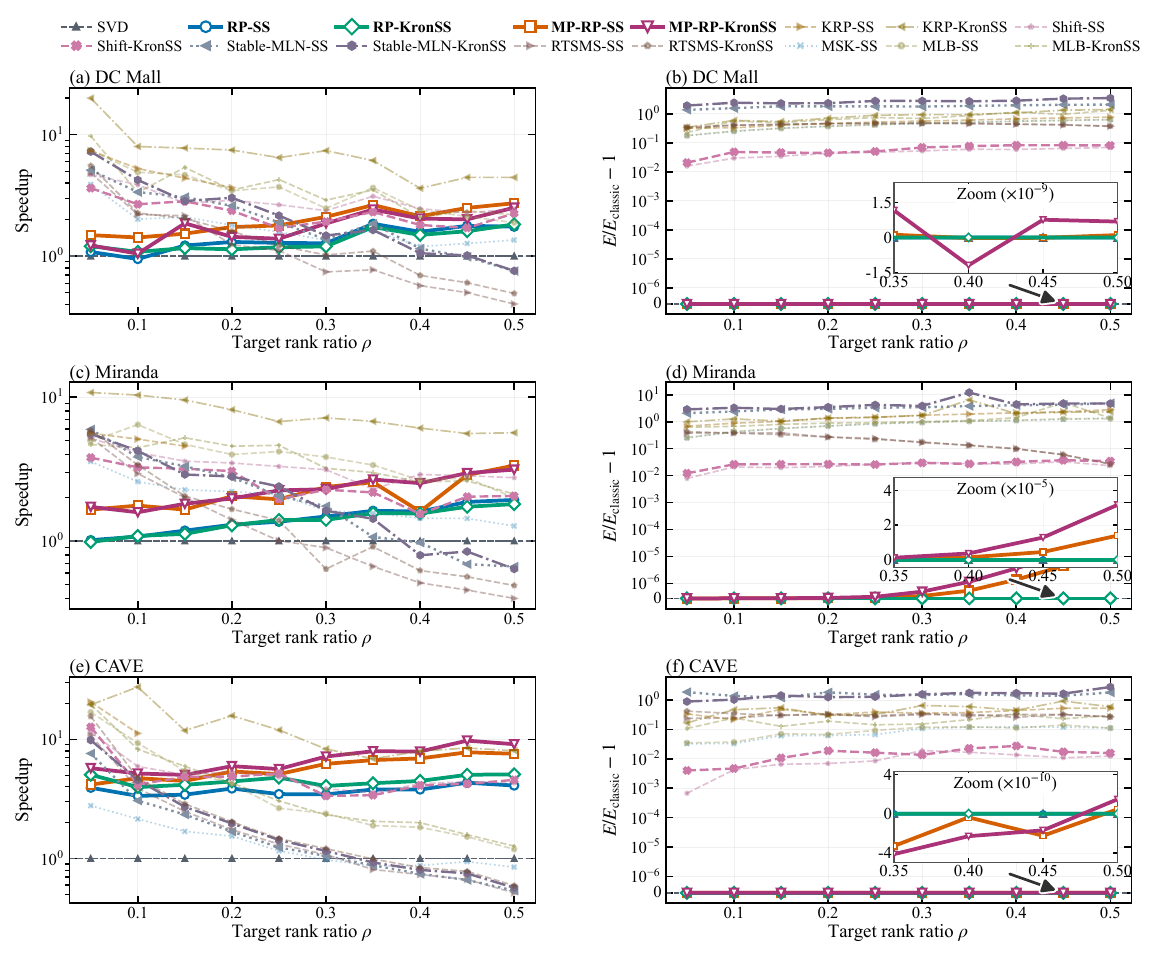}
 \caption{Real-data comparison: STHOSVD tests on DC Mall, Miranda and CAVE.}
 \label{fig:real-st}
\end{figure}

Fig.~\ref{fig:real-h} and~\ref{fig:real-st} report the results for DC Mall, Miranda and CAVE with target rank
$r_j=\max\{1,\operatorname{round}(\rho n_j)\}$, where
  $\rho=0.05,0.10,\ldots,0.50$.

\paragraph{Accuracy.} RP again matches the classical recovery error to numerical precision. MP-RP also remains close to the classical reference: the relative discrepancies on DC Mall and CAVE are below $0.000002\%$.  On Miranda, small deviations become visible at larger target rank ratios, but the largest discrepancy is only about $0.00317\%$. Thus the nearly classical recovery accuracy observed on synthetic tensors also carries over to these real data. We also observe that   occasional negative values of $E/E_{\mathrm{classic}}-1$ in the insets of Figs.~\ref{fig:real-h} and~\ref{fig:real-st}, which indicate   errors slightly below the classical ones. This is because classical HOSVD and STHOSVD provide quasi-optimal approximations, so their   errors need not be minimal at the prescribed ranks.

\paragraph{Time.} The proposed methods also accelerate the classical decompositions on real data. In the HOSVD tests, MP-RP-KronSS achieves speedups of about $5.5$--$6.3$ on DC Mall, $3.9$--$6.3$ on Miranda, and $20.1$--$21.9$ on CAVE. The corresponding STHOSVD speedups are $1.1$--$2.5$, $1.6$--$3.1$, and $5.0$--$9.8$. CAVE shows particularly clear additional gains from KronSS and mixed precision. The gains are more modest on DC Mall and Miranda in the STHOSVD tests, where   RP is sometimes close to the classical running time.

\clearpage

\section{Conclusions}\label{sec:conclusions}

In this work, we developed randomized preconditioning methods for HOSVD and STHOSVD, exploiting tensor unfoldings whose mode product dimensions can be enormous while their mode dimensions remain manageable. For each mode, we construct a preconditioner from a sketch, using numerical truncation to avoid inverting very small singular values, and then recover the dominant singular subspace from the unfolding itself. The large matrix products in the computation can be carried out in lower precision for further acceleration. For both decompositions, numerical truncation adds a term to the classical tail-energy bound that is quadratic in both the truncation threshold and the embedding distortion. This dependence is sharp in general. We also quantified the additional error due to the low-precision products. In tests on synthetic and real tensors, the proposed methods achieved speedups ranging from a fewfold to several dozenfold over classical HOSVD and STHOSVD. RP matched the classical reconstruction errors to numerical precision, while MP-RP achieved nearly the same accuracy.

\clearpage
\appendix
\begingroup
\small
\section{Rounding errors in Algorithm~\ref{alg:rp}}\label{app:rp-roundoff}
We follow the rounding model and matrix multiplication bound of Sect.~\ref{sec:mixed-rp},   with all operations at unit roundoff $u$. Write $\gamma_k=ku/(1-ku)$. A product with inner dimension $k$ and $ku<1$ satisfies $|\operatorname{fl}(AB)-AB|\le\gamma_k|A||B|$. Scalar division and Cholesky factorization obey the standard bounds~\cite[Chapters~3 and~10]{Higham2002}. For the backward-stable SVDs of $Y$ and $B_\epsilon$ in the algorithm, we represent the computed singular vector matrices by exactly orthogonal matrices differing by $O(u)$~\cite[\S4.9.1]{LAPACK1999} and absorb these corrections into the error terms.

\begin{proposition}\label{prop:rp-backward-stability}
Let $H\ne0$, suppose that~\eqref{eq:ose-standard} holds with $\mathscr U=\range(H)$, and take $\epsilon=cu\norm{Y}_2$ for the computed sketch $Y$. Let $U_B,\Sigma_B,V_B$ be the compact SVD factors of the computed $B_\epsilon$. Then there exist constants $c_0,C_0>0$ such that,
whenever $c\ge c_0$ and $C_0u\le1$, there exist
an orthonormal $U_H$   and a perturbation
$\Delta H$ satisfying
\begin{equation}\label{eq:rp-svd-backward}
 H+\Delta H=U_H\Sigma_BV_B\T,
 \qquad \norm{\Delta H}_F\le Cu\norm{H}_F.
\end{equation}
 $c_0,C_0$ depend only on the matrix dimensions, $\eta$, and $\norm{S}_2$; $C$ may also depend on the fixed $c$. All are independent of $\kappa_2(H)$. Therefore, Algorithm~\ref{alg:rp}  has the same order of normwise backward error as a direct backward stable SVD of $H$.

\end{proposition}

\begin{proof}
Write the computed sketch as $Y=SH+E_{\mathrm{sk}}$ and its backward stable SVD as $Y+E_{\mathrm{svd}}=U\Sigma V\T$. Write $E_s=E_{\mathrm{sk}}+E_{\mathrm{svd}}$; then $\bar Y:=SH+E_s=U\Sigma V\T$. The standard sketch and SVD error bounds, with $\norm{Y}_F=O(\norm{H}_F)$, imply $\norm{E_s}_F\le au\norm{H}_F$, where $a$ depends only on the dimensions and $\norm{S}_2$, independently of $\kappa_2(H)$. For the   SVD factors $V_\epsilon,\Sigma_\epsilon$ of $\bar Y$, write $P_\epsilon=V_\epsilon V_\epsilon\T$, $X_\epsilon=HV_\epsilon\Sigma_\epsilon^{-1}$, and let $\widetilde X_\epsilon$ be the computed preconditioned matrix:
\begin{equation}\label{eq:rp-scaled-roundoff}
 \widetilde X_\epsilon=X_\epsilon+E_x,\qquad
 \norm{E_x\Sigma_\epsilon}_F\le b u\norm{H}_F;
\end{equation}
to obtain~\eqref{eq:rp-scaled-roundoff}, the standard
multiplication and division bounds give
$|E_x\Sigma_\epsilon|\le\gamma_{n+1}|H||V_\epsilon|$, which implies $
\|E_x\Sigma_\epsilon\|_F
\le\gamma_{n+1}\|H\|_F\|V_\epsilon\|_F
=\gamma_{n+1}\sqrt t\,\|H\|_F
\le bu\|H\|_F
$ for some   $b$ depending only on the dimensions and independent  of $\kappa_2(H)$.

Next, the truncation rule gives $\sigma_t(\bar Y)>\epsilon=cu\norm{Y}_2$. The embedding and sketch error bounds yield $\norm{Y}_2\ge\alpha\norm{H}_F$ for some $\alpha>0$ independent of $\kappa_2(H)$, hence $\sigma_t(\bar Y)>c\alpha u\norm{H}_F$. Combining this with the bounds on $\norm{E_s}_F$ and $\norm{E_x\Sigma_\epsilon}_F$ gives

\begin{align}\label{eq:rp-error-bounds}
 \norm{E_sV_\epsilon\Sigma_\epsilon^{-1}}_2\le\frac{a}{c\alpha},
 \qquad \norm{E_x}_2\le\frac{b}{c\alpha}.
\end{align}
The identities $\bar YV_\epsilon\Sigma_\epsilon^{-1}=U_\epsilon$
and $SH+E_s=\bar Y$, together with the definition of
$X_\epsilon$, give $SX_\epsilon=U_\epsilon-E_sV_\epsilon\Sigma_\epsilon^{-1}$.
By~\eqref{eq:rp-error-bounds}, any unit vector $z$ satisfies
$1-a/(c\alpha)\le\norm{SX_\epsilon z}_2
\le\sqrt{1+\eta}\norm{X_\epsilon z}_2$,
where the second inequality follows from the embedding
bound on $\range(H)$. Since $\widetilde X_\epsilon=X_\epsilon+E_x$,
it thus follows that

\[
 \sigma_{\min}(\widetilde X_\epsilon)
 \ge\frac{1-a/(c\alpha)}{\sqrt{1+\eta}}-\frac{b}{c\alpha}.
\]
 Therefore, for sufficiently large $c_0$ and $c\ge c_0$,
$\sigma_{\min}(\widetilde X_\epsilon)$ has a positive constant
lower bound. Similarly, the other side of the embedding bound
gives $\norm{\widetilde X_\epsilon}_2=O(1)$.
Both bounds are independent of $\kappa_2(H)$, so
$\widetilde X_\epsilon\T\widetilde X_\epsilon$ is bounded
and uniformly positive definite.
The standard Gram and Cholesky error bounds then give
$R_\epsilon\T R_\epsilon
=\widetilde X_\epsilon\T\widetilde X_\epsilon+F$,
with $\norm{F}_2=O(u)$ and $\norm{R_\epsilon^{-1}}_2=O(1)$.
It follows that
$(\widetilde X_\epsilon R_\epsilon^{-1})\T
(\widetilde X_\epsilon R_\epsilon^{-1})=I+O(u)$,
and its polar factor $Q$ has orthonormal columns with
$\norm{\widetilde X_\epsilon R_\epsilon^{-1}-Q}_2=O(u)$.
Accounting for the multiplication error $E_b$ from forming $B_\epsilon$, we write
$B_\epsilon=R_\epsilon\Sigma_\epsilon V_\epsilon\T+E_b$.
By~\eqref{eq:rp-scaled-roundoff} and
$HP_\epsilon=X_\epsilon\Sigma_\epsilon V_\epsilon\T$,
\begin{align*}
 HP_\epsilon-QB_\epsilon
 &=(\widetilde X_\epsilon-E_x)\Sigma_\epsilon V_\epsilon\T
   -Q(R_\epsilon\Sigma_\epsilon V_\epsilon\T+E_b)\\
 &=-E_x\Sigma_\epsilon V_\epsilon\T
   +(\widetilde X_\epsilon-QR_\epsilon)\Sigma_\epsilon V_\epsilon\T-QE_b\\
 &=-E_x\Sigma_\epsilon V_\epsilon\T
 +(\widetilde X_\epsilon R_\epsilon^{-1}-Q)R_\epsilon\Sigma_\epsilon V_\epsilon\T-QE_b;
\end{align*}
each term on the right has Frobenius norm $O(u)\norm{H}_F$.
Indeed, the bound for $E_x\Sigma_\epsilon V_\epsilon\T$ follows
from~\eqref{eq:rp-scaled-roundoff} and the orthonormal columns
of $V_\epsilon$. For the second term, use
$\norm{\widetilde X_\epsilon R_\epsilon^{-1}-Q}_2=O(u)$,
$\norm{R_\epsilon}_2^2=\norm{\widetilde X_\epsilon}_2^2+O(u)=O(1)$,
and $\norm{\Sigma_\epsilon}_F\le\norm{\bar Y}_F=O(\norm{H}_F)$.
For the third term,
$\norm{QE_b}_F=\norm{E_b}_F
\le Cu\norm{R_\epsilon}_2\norm{\Sigma_\epsilon}_F
\norm{V_\epsilon\T}_2=O(u)\norm{H}_F$.

It remains to bound $\norm{H(I-P_\epsilon)}_F$.
The embedding condition and $SH=\bar Y-E_s$ give
\[
 \begin{aligned}
 \sqrt{1-\eta}\,\norm{H(I-P_\epsilon)}_F
 &\le\norm{SH(I-P_\epsilon)}_F
 \le\norm{\bar Y(I-P_\epsilon)}_F+\norm{E_s}_F\\
 &\le\sqrt{n-t}\,\epsilon+\norm{E_s}_F;
 \end{aligned}
\]
the last step uses the fact that all $n-t$ discarded singular values are at most $\epsilon$. Substituting $\epsilon=cu\norm{Y}_2=O(u)\norm{H}_F$ and $\norm{E_s}_F\le au\norm{H}_F$ yields $\norm{H(I-P_\epsilon)}_F=O(u)\norm{H}_F$. Altogether, $\norm{H-QB_\epsilon}_F\le\norm{H(I-P_\epsilon)}_F+\norm{HP_\epsilon-QB_\epsilon}_F=O(u)\norm{H}_F$. Finally, standard backward stability of the SVD of $B_\epsilon$ gives $B_\epsilon+E_B=U_B\Sigma_BV_B\T$, with $\norm{E_B}_F=O(u)\norm{B_\epsilon}_F$. Taking $U_H=QU_B$ and $\Delta H=Q(B_\epsilon+E_B)-H$ proves the result.
\end{proof}
\endgroup

\begingroup
\footnotesize
\noindent\textbf{Use of AI tools.}
The research ideas, algorithm design, theoretical objectives, experimental design, and manuscript structure were developed by the authors. ChatGPT and Codex assisted in an author-led iterative process to examine and refine proof arguments, write code for numerical experiments and Lean formalizations, and improve the clarity  of the manuscript at the sentence level. The authors critically reviewed the AI-assisted material, made all final decisions, and take full responsibility for the content of the paper.\par
\endgroup

\providecommand{\doi}[1]{\url{https://doi.org/#1}}

\end{document}